\documentclass[letterpaper, 12pt, oneside]{report}

    \usepackage[left=1.5in, right=1in, top=1in, bottom=1in, dvips, pdftex]{geometry}
    
    \usepackage{fancyhdr} 
    \usepackage{amsmath,amsthm,amsfonts,amssymb,amscd,enumerate,graphicx,etex,diagrams,ifsym}
    \usepackage{setspace} 
    \usepackage[linktocpage,bookmarksopen,bookmarksnumbered,
                pdftitle={My Doctoral Thesis},
                pdfauthor={Me},%
                pdfsubject={UWM Ph.D. Doctoral Thesis},%
                pdfkeywords={UWM Ph.D. Doctoral Thesis}]{hyperref}
    
\usepackage{paralist} 

     \usepackage[toc,title]{appendix} 
     \usepackage{MnSymbol}
     \usepackage{enumerate}
     \usepackage{amsfonts}

\theoremstyle{plain}

\newtheorem{Theorem}[subsection]{Theorem} 
\newtheorem{Proposition}[subsection]{Proposition}
\newtheorem{Corollary}[subsection]{Corollary}
\newtheorem{Lemma}[subsection]{Lemma}
\newtheorem{Claim}[subsection]{Claim}
\newtheorem{Definition}[subsection]{Definition}
\newtheorem{Remark}[subsection]{Remark}
\newtheorem{Example}[subsection]{Example}
\newtheorem{Fact}[subsection]{Fact}
\newtheorem{Problem}[subsection]{Problem}

\newenvironment{Sketch}{\par\noindent{\sc Sketch of Proof}\quad}{\hfill\qed\par\smallskip}
\newenvironment{Proof}{\par\noindent{\sc Proof}\quad}{\hfill\qed\par\smallskip}

\newenvironment{RestateTheorem}[3]{\par\vspace{12pt}\noindent{\bf #1~\ref{#2}}(#3){\bf .}\it}{\par\vspace{12pt}}

\numberwithin{equation}{section}
\numberwithin{figure}{section}
\newcommand{\n}{\vspace{12pt}} 
      
\newcommand{\newchapter}[3] 
	{                           
        \chapter[#2]{#3}
        \chaptermark{#1}
        \thispagestyle{myheadings}
	}

\newcommand{\newappendix}[3] 
	{                           
                
        \chapter[#2]{#3}
        \chaptermark{#1}
        \thispagestyle{myheadings}
	}

\newcommand{\BN}{\mathbb{N}}
\newcommand{\BZ}{\mathbb{Z}}

\newcommand{\BS}{\mathbb{S}}

\newcommand{\BI}{\mathbb{I}}

\newcommand{\rto}{\rightarrow}
\newcommand{\lto}{\leftarrow}

\fancypagestyle{vita}{
\fancyhf{}
\fancyhead[R]{\thepage}

}

        {\begin{enumerate}[#1]}{\end{enumerate}%
         \vspace{-.6\baselineskip}}
        {\begin{compactenum}[#1]}{\end{compactenum}}

\fancypagestyle{references}{
\fancyhf{}
\fancyhead[R]{\thepage}

}

\numberwithin{equation}{section}

\begin{document}
    

    \pagenumbering{roman}
    \pagestyle{plain}

    %
    %
\thispagestyle{empty}
    \singlespacing

    ~\vspace{1.5in} 
    \begin{center}

        \begin{huge}
            SOME RESULTS ON PSEUDO-COLLAR STRUCTURES ON HIGH-DIMENSIONAL MANIFOLDS
        \end{huge}\\\n
        By\\\n
        {\sc Jeffrey Joseph Rolland}\\
        A Dissertation Submitted in\\\n
        Partial Satisfaction of the\\\n
        Requirements for the Degree of\\\n\n
        DOCTOR OF PHILOSOPHY\\\n
        in\\\n
        MATHEMATICS\\\n\n
        at the \\\n
        UNIVERSITY OF WISCONSIN-MILWAUKEE\\\n
                
        May, 2015\\
	
	\end{center}
\newpage

%

    
%
%
%
%
%
%
%
%
%

    \newpage
       
    %
    %
    %
    
    %

    
    \chapter*{\vspace{-1.5in}Abstract}
    
    \begin{center}
    SOME RESULTS ON PSEUDO-COLLAR STRUCTURES ON HIGH-DIMENSIONAL MANIFOLDS\\\n
      
      by\\\n
      
      Jeffrey Joseph Rolland\\\n
      
      The University of Wisconsin-Milwaukee, 2015\\
      Under the Supervision of Prof. Craig R. Guilbault
    \end{center}
    
    \doublespacing
    In this paper, we provide expositions of	 Quillen's plus construction for high-dimensional smooth manifolds and the solution to the group extension problem. We then develop a geometric procedure due for producing a ``reverse'' to the plus construction, a construction called a \textit{semi-s-cobordism}. We use this reverse to the plus construction to produce ends of manifolds called \textit{pseudo-collars}, which are stackings of semi-h-cobordisms. We then display a technique for producing ``nice'' one-ended open manifolds which satisfy two of the necessary and sufficient conditions for being pseudo-collarable, but not the third. Finally, we recall a different, but very difficult to enact in practice, procedure due to J.-C. Hausmann and P. Vogel which enumerates the class of all semi-s-cobordisms for a given closed manifold, but does not tell when this set is non-empty. We show a connection between Hausmann-Vogel's technique for producing semi-s-cobordisms and our technique for producing semi-s-coborisms.

    \n\n\n\n
    
    \newpage
    \singlespacing 
    %
    %
    %
        

    ~\\[3.625in] 
    \begin{center}
      \copyright ~by Jeffrey Rolland, 2015.\\
      All rights reserved.
    \end{center}
%
    \newpage
%
%
    %
    %
    
    ~\\[3.625in] 
    \centerline{For Rev\'e, and for my parents, Thomas and Mary Lou}
    
    \newpage

    %
    %
    
    \tableofcontents
    
    \newpage

    %

%

%



%
   
    %
    %

    \chapter*{\vspace{-1.5in}Acknowledgements}

        The author would like to acknowledge the immense input and direction of his major professor, Craig Guilbault. He would also like to acknowledge the help and assistance (roughly in order of assistance) of Jeb Willenbring of the University of Wisconsin - Milwaukee, Jack Schmidt of the University of Kentucky, ``clam'' of the DALnet IRC channel \#macintosh, Ric Ancel of the University of Wisconsin - Milwaukee, Boris Okun of the University of Wisconsin - Milwaukee, Chris Hruska of the University of Wisconsin - Milwaukee, Alejandro Adem of the University of British Columbia (Canada), Ross Geoghegan of Binghamton Univeristy, Mike Hill of the University of Virginia, Allen Bell of the University of Wisconsin - Milwaukee, Derek Holt of the University of Warwick (UK), Marston Conder of the University of Auckland (New Zealand), Jason Manning of the University at Buffalo, and Brita Nucinkis of the University of London.
    
    \newpage

    %
    %

    \pagestyle{fancy}
    \pagenumbering{arabic}
    \doublespacing
       
    %
    %

    \newchapter{INTRODUCTION}{INTRODUCTION}{INTRODUCTION}
    \label{ch:1:Intro}
    
    	
	
	\section[Informal Overview and Historical Motivation]{Informal Overview and Historical Motivation}
    \label{sec:Intro:Motivation}
    We work in the category of smooth manifolds, but all our results apply equally well to the categories of PL and topological manifolds. The manifold version of Quillen's plus construction provides a way of taking a closed  smooth manifold $M$ of dimension $n \ge 5$ whose fundamental group $G = \pi_1(M)$ contains a perfect normal subgroup $P$ which is the normal closure of a finite number of elements and producing a compact cobordism $(W,M,M^+)$ to a manifold $M^+$ whose fundamental group is isomorphic to $Q = G/P$ and for which $M^+ \hookrightarrow W$ is a simple homotopy equivalence. By duality, the map $f:M \rightarrow M^+$ given by including $M$ into $W$ and then retracting onto $M^+$ induces an isomorphism $f_*:H_*(M;\mathbb{Z}Q) \rightarrow H_*(M^+;\mathbb{Z}Q)$ of homology with twisted coefficients. By a clever application of the s-Cobordism Theorem, such a cobordism is uniquely determined by $M$ and $P$ (see \cite{Freedman-Quinn} P. 197).

\indexspace

In ``Manifolds with Non-stable Fundamental Group at Infinity I'' \cite{Guilbault1}, Craig Guilbault outlines a structure to put on the ends of an open smooth manifold $N$ with finitely many ends called a \textit{pseudo-collar}, which generalizes the notion of a collar on the end of a manifold introduced in Siebenmann's dissertation \cite{Siebenmann}. A pseudo-collar is defined as follows. Recall that a manifold $U^n$ with compact boundary is an open collar if $U^n \approx \partial U^n \times [0,\infty)$; it is a homotopy collar if the inclusion $\partial U^n \hookrightarrow U^n$ is a homotopy equivalence. If $U^n$ is a homotopy collar which contains arbitrarily small homotopy collar neighborhoods of infinity, then we call $U^n$ a \textit{pseudo-collar}. We say that an open $n$-manifold $N^n$ is collarable if it contains an open collar neighborhood of infinity, and that $N^n$ is \textit{pseudo-collarable} if it contains a pseudo-collar neighborhood of infinity.

\indexspace

Each pseudo-collar admits a natural decomposition as a sequence of compact cobordisms $(W,M,M_-)$, where $W$ is a semi-h-cobordism (see Definition \ref{defsemi-s-cob} below). It follows that the cobordism $(W,M_-,M)$ is a one-sided h-cobordism (a plus cobordism if the homotopy equivalence is simple). (This somewhat justifies the use of the symbol ``$M_-$'' for the right-hand boundary of a semi-h-cobordism, a play on the traditional use of $M^+$ for the right-hand boundary of a plus cobordism.)

\indexspace

The general problem of a reverse to Quillen's plus construction in the high-dimensional manifold category is as follows.

\begin{Problem}[Reverse Plus Problem]
Suppose $G$ and $Q$ are finitely-presented groups and $\Phi: G \twoheadrightarrow Q$ is an onto homomorphism with $\ker(\Phi)$ perfect. Let $M^n$ ($n \ge 5$) be a closed smooth manifold with $\pi_1(M) \cong Q$.

\indexspace

Does there exist a compact cobordism $(W^{n+1}, M, M_-)$ with 

\begin{diagram}[size=14.5pt]
1 & \rTo & \ker(\iota_{\#}) & \rTo & \pi_1(M_-) & \rTo^{\iota_{\#}} & \pi_1(W) & \rTo & 1
\end{diagram}

equivalent to

\begin{diagram}[size=14.5pt]
1 & \rTo & \ker(\Phi) & \rTo & G & \rTo^{\Phi} & Q & \rTo & 1
\end{diagram}

and $M \hookrightarrow W$ a (simple) homotopy equivalence. 
\end{Problem}

Notes:
\begin{itemize}

\item The fact that $G$ and $Q$ are finitely presented forces $\ker(\Phi)$ to be the normal closure of a finite number of elements. (See, for instance, \cite{Guilbault1} or \cite{Siebenmann}.)

\item Closed manifolds $M^n$ ($n \ge 5$) in the various categories with $\pi_1(M)$ isomorphic to a given fintely presented group $Q$ always exist. In the PL category, one can simply take a presentation 2-complex for $Q$, $K$, embed $K$ in $\BS^{n+1}$, take a regular neighborhood $N$ of $K$ in $\BS^{n+1}$, and let $M = \partial N$. Similar procedures exist in the other categories.

\end{itemize}

The following terminology was first introduced in \cite{Hausmann1}.

\begin{Definition} [1-Sided e-Cobordism] \label{defsemi-s-cob}
Let $N^n$ be a compact smooth manifold. A \textbf{1-sided e-cobordism} $(W,N,N_-)$ is a cobordism so that $N \hookrightarrow W$ is a homotopy equivalence (necessarily simple if $e = s$ and not necessarily simple if $e = h$). [A 1-sided e-cobordism  $(W,N,N_-)$ is so-named presumably because it is ``half an e-cobordism''].
\end{Definition}

One wants to know under what circumstances 1-sided e-cobordisms exists, and, if they exists, how many there are.

\indexspace

There are some cases in which 1-sided s-cobordisms are known not to exist. For instance, if $P$ is fintely presented and perfect but not superperfect, $Q = \langle e \rangle$, and $M = \BS^n$, then a solution to the Reverse Plus Problem would produce an $M_-$ that is a homology sphere. But it is a standard fact that a manifold homology sphere must have a a superperfect fundamental group! (See, for instance, \cite{Kervaire}.) (The definition of superperfect will be given in Definition \ref{defsuperperfect}.)

\indexspace

The key point is that the solvability fo the Reverse Plus Problem depends not just upon the group data, but also upon the manifold $M$ with which one begins. For instance, one could start with a group $P$ which is finitely presented and perfect but not superperfect, let $N_-$ be a manifold obtained from the boundary of a regular neighborhood of the embedding of a presentation 2-complex for $P$ in $\BS^{n+1}$, and let $(W, N_-, N)$ be the result of applying Quillen's plus construction to to $N_-$ with respect to all of $P$. Then again $Q = \langle e \rangle$ and $\Phi: P \twoheadrightarrow Q$ but $N$ clearly admits a semi-s-cobordism, namely $(W, N, N_-)$ (however, of course, we cannot have $N$ a sphere or $N_-$ a homology sphere).

\indexspace

Hausmann and Vogel's work in \cite{Hausmann1}, \cite{Hausmann2}, and \cite{H-V} provides a partial solution to the Reverse Plus Problem in the case the kernel group is \textit{locally perfect}, that is, when every element of the kernel group is contained in a finitely generated perfect subgroup. They set up an obstruction theory which puts solutions to a given Reverse Plus Problem in one-to-one correspondence with a carefully defined collection of maps $\{X_M, BG^+\}$. So, the Reverse Plus Problem asks whether this set is non-empty for a specific set of initial data. (As noted above, the set $\{X_M, BG^+\}$ may well be empty.) Our Theorem \ref{thmsemi-s-cob} bypasses their theory and gives a direct method for constructing a solution to the Reverse Plus Problem in certain stiuations. Their theory was unknown to us at the time we proved Theorem \ref{thmsemi-s-cob}, but, in retrospect, our result can be viewed as a proof that their classifying set is non-empty in certain situations - an issue which they addressed in only a few select cases.

\indexspace

Per \cite{D-T2}, the group $G = \langle x, y, t\ |\ y = [y,y^x], txt^{-1} = y \rangle$ admits a subgroup $H = \langle x, y\ |\ y = [y,y^x] \rangle$ with the properties that the normal closure of $H$ in $G$, $ncl(H;G)$, is perfect, $G/ncl(H;G) \cong \BZ$, and that no finitely generated subgroup of $ncl(H;G)$ contains a non-trivial perfect subgroup. So, $ncl(H;G)$ is perfect but not locally perfect. Thus, there are instances of the Reverse Plus Problem to which Hausmann and Vogel's work does not apply at all. Unfortunatly, Theorem \ref{thmsemi-s-cob} does not apply to these instances, either.
	
	 \section[Main Results]{Statements of the Main Results}
     \label{sec:LabelForChapter1:Section2}
     \begin{Theorem}[Existence of 1-sided s-cobordisms] \label{thmsemi-s-cob}
Given $1 \rightarrow S \rightarrow G \rightarrow Q \rightarrow 1$ where $S$ is a finitely presented superperfect group, $G$ is a semi-direct product of $Q$ by $S$, and any $n$-manifold $N$ with $n \ge 6$ and $\pi_1(M) \cong Q$, there exists a solution $(W, N, N_-)$ to the Reverse Plus Problem for which $N \hookrightarrow W$ is a simple homotopy equivalence.
\end{Theorem}

One of the primary motivations for Theorem \ref{thmsemi-s-cob} is that it provides a ``machine'' for constructing interesting pseudo-collars. As an application, we use it to prove:

\begin{Theorem}[Uncountably Many Pseudo-Collars on Closed Manifolds with the Same Boundary and Similar Pro-$\pi_1$]  \label{thmpseudo-collars}
Let $M^n$ be a closed smooth manifold ($n \ge 6$) with $\pi_1(M) \cong \BZ$ and let $S$ be the fintely presented group $V*V$, which is the free preduct of 2 copies of Thompson's group $V$. Then there exists an uncountable collection of pseudo-collars $\{N^{n+1}_{\omega}\ |\ \omega \in \Omega\}$, no two of which are homeomorphic at infinity, and each of which begins with  $\partial N^{n+1}_{\omega} = M^n$ and is obtained by blowing up countably many times by the same group $S$. In particular, each has fundamental group at infinity that may be represented by an inverse sequence

\begin{diagram}[size=14.5pt]
\BZ & \lOnto^{\alpha_1} & G_{1} & \lOnto^{\alpha_2} & G_{2} & \lOnto^{\alpha_3} & G_{3} & \lOnto^{\alpha_4} & \ldots \\
\end{diagram}

with $ker(\alpha_i) = S$ for all $i$.
\end{Theorem}

As noted above, our work allows the construction of a wide variety of pseudo-collars. In a related but different direction, we expand upon an example by Guilbault and Tinsley found in \cite{G-T2}, by describing a procedure for constructing a wide variety of very nice ends which nevertheless do not admit pseudo-collar structures. More specifically, we prove:

\begin{Theorem}[Existence of Non-Pseudo-Collarable ``Nice'' Manifolds] \label{thm-non-pcm}
Let $M^n$ be an orientable, closed manifold ($n \ge 6$) such that $\pi_1(M)$ contains an element $t_0$ of infinite order and $\pi_1(M)$ is hypo-Abelian (defined in Section 5.1). Then there exists a 1-ended, orientable manifold $W^{n+1}$ with $\partial W = M$ in which all clean neighborhoods of infinity have finite homotopy type, but which does not have perfectly semistable fundamental group at infinity. Thus, $W^{n+1}$ is absolutely inward tame but not pseudocollable.
\end{Theorem}

	 \section[Notation]{Some Notational Conventions}
     \label{sec:LabelForChapter1:Section3}
     Throughout, $\mathbb{N}$ denotes the set of natural numbers (note $0 \in \mathbb{N}$), $\mathbb{Z}$ denotes the integers, $\mathbb{Q}$ denotes the rational numbers, $\mathbb{R}$ denotes the real numbers, and $\mathbb{C}$ denotes the complex numbers.

\indexspace

$\mathbb{I}$ is the closed unit interval [0,1]. $\mathbb{D}^n = \{{\bf x} \in \mathbb{R}^n\ |\ ||{\bf x}|| \le 1\}$ is the n-disk. $\mathbb{S}^n = \partial\mathbb{D}^{n+1}$ is the n-sphere.

\indexspace

$\mathbb{Z}_n$ denotes the integers mod n. $\mathbb{F}_{p^n}$ denotes the field with $p^n$ elements, $p$ a prime.

        %
    
        
        
            
        
        

    %
    %

    \newchapter{PLUS CONSTRUCTION AND GROUP EXTENSION PROBLEM}{THE PLUS CONSTRUCTION AND THE GROUP EXTENSION PROBLEM}{THE PLUS CONSTRUCTION FOR HIGH DIMENSIONAL MANIFOLDS AND THE SOLUTION TO THE GROUP EXTENSION PROBLEM}
    \label{sec:LabelForChapter2}

			In this chapter, we present an overview of the plus construction and the solution to the group extension problem. Neither of these constructions are original work due to the author; they are, however, fundamental background material necessary for understanding the author's work. 1-sided s-cobordisms have the property that, when read in one direction (sometimes called ``semi-s-cobordisms'') are ``reverse plus constructions'', but when read in the other direction are ``forward'' plus constructions, so we provide some information on 1-sided s-cobordisms, as well as 1-sided h-cobordisms. Also, solving the group extension problem, as well as some of the subtleties involved in solving the group extension problem (such as the distinction between a semi-direct product and other solutions utilizing the same abstract kernel), are crucial to understanding the author's work on a reverse to the plus construction, so some information on the solution to the group extention problem is included. The reader familiar with both of these techniques may safely skip this chapter.
			    
        \section[Review of The Plus Construction]{Review of The Plus Construction in High Dimensions}
        \label{sec:LabelForChapter2:Section1}
        In this section, we give an exposition of the Manifold Plus Construction in High Dimensions. The Plus Construction is a well-known and important work originally due to  Quillen in \cite{Quillen}. The manifold version is a bit more complicated, due to framing issues, and seems to be a part of the folklore: see \cite{G-T5}.

\begin{Theorem}[The Manifold Plus Construction in High-Dimensions]
Let $M^n$ be a closed smooth manifold of dimension $n \ge 5$ such that $G = \pi_1(M,\star)$ contains a perfect normal subgroup $P$ which is normally generated in $G$ by a finite number of elements. Let $\Phi: G \rightarrow Q = G/P$ be the quotient map. Then there is a compact cobordism $(W,M,M^+)$ to a manifold $M^+$ with the following properties: 
\begin{itemize}
\item the short exact sequence 

\begin{diagram}[size=14.5pt]
1 & \rTo & \ker(\iota_{\#}) & \rTo & \pi_1(M_-) & \rTo^{\iota_{\#}} & \pi_1(W) & \rTo & 1
\end{diagram}

is equivalent to

\begin{diagram}[size=14.5pt]
1 & \rTo & P & \rTo & G & \rTo^{\Phi} & Q & \rTo & 1
\end{diagram}
\item the inclusion $\iota: M^+ \hookrightarrow W$ is a simple homotopy equivalence
\item the map $f: M \rightarrow M^+$ given by including $M$ into $W$ and then retracting $W$ onto $M^+$ induces an isomorphism $f_*:H_*(M;\mathbb{Z}Q) \rightarrow H_*(M^+;\mathbb{Z}Q)$ of homology with twisted coefficients
\item The manifold $M^+$ is unique up to a diffeomorphism, and $W$ is unique up to a diffeomorphism rel $M$ and $M^+$.
\end{itemize}
\end{Theorem}

Next, we turn our attention to the functoriality of the plus construction. Our basic approach is based upon \cite{Loday}, adopted for the manifold categories.

\begin{Theorem} \label{thmloday-1}
Let $M^n$ be a closed smooth manifold of dimension $n \ge 5$ such that $G = \pi_1(M,\star)$ contains a perfect normal subgroup $P$ which is normally generated in $G$ by a finite number of elements. Suppose $N^n$ is a closed smooth manifold and there is a smooth map $f: M \rightarrow N$ such that $f_{\#}(P) =\ \langle e \rangle\ \le \pi_1(N)$. Then there are smooth maps $\iota_M: M \rightarrow M^+$ and $f^+: M^+ \rightarrow N$ such that $f = f^+ \circ \iota_M$
\end{Theorem}

\begin{Corollary}[The Plus Construction is Functorial]
Given smooth manifolds $M$ and $N$ and a smooth map $f:M \rightarrow N$, suppose there are perfect normal subgroups $P$ of $G = \pi_1(M)$ and $P'$ of $G' = \pi_1(N)$ such that $P$ and $P'$ are normally finitely generated in $G$ and $G'$ respectively and $f_{\#}(P) \le P'$. Then there are smooth maps $\iota_M: M \rightarrow M^+$, $\iota_N: N \rightarrow N^+$, and $f^+:M^+ \rightarrow N^+$ such that $f^+ \circ \iota_M = \iota_N \circ f$, where $M^+$ is the result of the plus construction applied to $M$ with respect to $P$ and $N^+$ the result of the plus construction applied to $N$ with respect to $P'$.
\end{Corollary}

\begin{Definition} \label{def1sided-cob}
A cobordism $(W,M,M^+)$ or $(W,M^+,M)$ is a \textbf{1-sided h-cobordism} if $M^+ \hookrightarrow W$ is a homotopy equivalanence. A cobordism $(W,M,M^+)$ or $(W,M^+,M)$ is a \textbf{1-sided s-cobordism} if $M^+ \hookrightarrow W$ is a simple homotopy equivalanence.
\end{Definition}

\begin{Definition}
Let $(W,M,M^+)$ or $(W,M^+,M)$ be a 1-sided s-cobordism. Then $(W,M^+,M)$ is called a \textbf{semi-s-cobordism} (although this term isn't used much in the modern literature) and $(W,M,M^+)$ is called a plus cobordism.
\end{Definition}

\indexspace

There is an analog of the Plus Construction in the CW complex category; see, for instance, Proposition 4.40 on page 374 in \cite{Hatcher} and its generalization in the paragraph immediately following the proof. We state this generalization for completeness.

\begin{Theorem}
Let $X$ be a connected CW complex with $P \le \pi(X)$ a perfect subgroup. Then there is a CW complex $X^+$ with $\pi(X^+) \cong \pi_1(X)/P$ and a map $f: X \rightarrow X^+$ inducing a quotient map on fundamental groups and isomorphisms on all homology groups with $\BZ Q$ coefficients.
\end{Theorem}
            
        \section[Review of The Group Extension Problem]{Review of The Solution to the Group Extension Problem}
        \label{sec:LabelForChapter2:Section2}
        In this section, we give an exposition of the Group Extension Problem. The following presentation of the solution to the Group Extension Problem is based upon \cite{MacLane} Chapter IV, section 8; the reader is referred there for more details. We include this presentation for completeness.

\begin{Definition}
Let $K$ and $Q$ be given groups. We say that a group $G$ \textit{solves the group extension problem for Q and K} or that \textit{G is an extension of Q by K} [WARNING: some authors use the reverse terminology and say G is an extension of K by Q] if there exists a short exact sequence 

\indexspace

\begin{diagram}
1 & \rTo & K & \rTo^\iota  & G & \rTo^\sigma  & Q & \rTo & 1 \\
\end{diagram}
\end{Definition}

Let $Aut(K)$ denote the automorphism group of $K$. Define $\mu: K \rightarrow Aut(K)$ to be $\mu(k)(k') = kk'k^{-1}$. Then the image of $\mu$ in $Aut(K)$ is called \textit{the inner automorphism group of K}, $Inn(K)$. The inner automorphism group of a group $K$ is always normal in $Aut(K)$. The quotient group $Aut(K)/Inn(K)$ is called the \textit{outer automorphism group} $Out(K)$. The kernel of $\mu$ is called the \textit{center of K}, $Z(K)$; it is the set of all $k \in K$ such that for all $k' \in K, kk'k^{-1} = k'$. One has the exact sequence

\indexspace

\begin{diagram}
1 & \rTo & Z(K) & \rTo  & K & \rTo^\mu  & Aut(K) & \rTo^\alpha & Out(K) & \rTo & 1 \\
\end{diagram}

\indexspace

Two group extensions $1 \rightarrow K \rightarrow G \rightarrow Q \rightarrow 1$ and $1 \rightarrow K \rightarrow G' \rightarrow Q \rightarrow 1$ are said to be \textit{congruent} if and only if there is an isomorphism $\gamma: G \rightarrow G'$ such that the following diagram commutes:

\indexspace

\begin{diagram}
1       & \rTo & K      & \rTo^\iota  & G           & \rTo^\sigma  & Q      & \rTo & 1 \\
\dTo^=  &      & \dTo^= &             & \dTo^\gamma &              & \dTo^= &      & \dTo^= \\
1       & \rTo & K      & \rTo^{\iota'} & G'          & \rTo^{\sigma'} & Q      & \rTo & 1 \\
\end{diagram}

\indexspace

Any group extension of $Q$ by $K$, $1 \rightarrow K \rightarrow G \rightarrow Q \rightarrow 1$, determines a homomorphism $\theta: G \rightarrow Aut(K)$ determined by conjugation: $\theta(g)(k') = \iota^{-1}(g\iota(k')g^{-1})$, if $\iota: K \rightarrow G$ is the inclusion map. Note this is well-defined, as $g\iota(k')g^{-1} \in \iota(K)$, as $\iota(K)$ is normal in $G$. Let $\xi: Aut(K) \rightarrow Out(K)$ denote the projection map. Note $\iota(K) \subseteq ker(\xi \circ \theta)$, as follows. First, observe $\theta(\iota(k))(k') = \iota^{-1}(\iota(k)\iota(k')\iota(k)^{-1}) = kk'k^{-1} = \mu(k)(k')$ meaning $\theta \circ \iota = \mu$. Next, see that, since $\mu(K) = Inn(K)$, we have $(\xi \circ \theta)(\iota(K)) = \xi((\theta \circ \iota)(k)) = \xi(\mu(K)) = \xi(Inn(K)) =\ <e>$. Therefore, as $ker(\sigma) = \iota(K)  \subseteq ker(\xi \circ \theta)$, we have a derived homomorphism $\psi: Q \rightarrow Out(K)$.

\indexspace

So, any group extension determines a homomorphism $\psi: Q \rightarrow Out(K)$. We call such a homomorphism an \textit{outer action of Q on K}. The homomorphism $\psi$ records the way in which $K$ appears as a normal subgroup of $G$. A pair of groups K and Q together with an outer action $\psi$ of Q on K is called an \textit{abstract kernel}.

\indexspace

The general problem of group extensions is to classify all group extensions up to congruence. Note that congruent extensions determine the same outer action.

\begin{Theorem}[Obstructions to Group Extensions]
Given a abstract kernel $(Q,K,\psi)$, interpret the center of $K$, $Z(K)$, as a $Q$-module, with the action $q.z = \phi(q)(z)$ for any choice of automorphism $\phi: K \rightarrow K$ with $\phi \cdot Inn(K) = \psi$. Then we may assign a cohomology class which vanishes if and only if $(Q,K,\psi)$ gives rise to a group extension.
\end{Theorem}

\begin{Theorem}[Classification of Group Extensions]
If an abstract kernel $(Q,K,\psi)$ has 0 obstruction, then the set of congruence classes of extensions with abstract kernel $(Q,K,\psi)$ is in bijective correspondence with the set $H^2(Q;Z(K))$, where $Z(K)$ has the module structure given in Theorem 2.3.7. This correspondence associates the $0 \in H^2(Q;Z(K))$ with the semi-direct product of $Q$ by $K$ with the given outer action (see below).
\end{Theorem}

It is important to note that this theorem only classifies group extensions up to congruence, not isomorphism. That is to say, two group extensions $G$ and $G'$ may use the same abstract kernel but different elements of $H^2(Q; Z(K))$ and therefore would not be congruent, but may still be isomorphic as groups. Indeed, it is theoretically possible that two group extensions may use different abstract kernels and still give rise to isomorphic extensions.

\begin{Definition} \label{defsemi-dir-prod}
A \textbf{semi-direct product} of $Q$ by $K$, $G = K \rtimes Q$, is a group extension of $Q$ by $K$ which splits, that is, for which there exist homomorphisms $j: G \rightarrow K$ and $k: Q \rightarrow G$ such that $\iota \circ j = id_K$ and $k \circ \sigma = id_Q$. 
\end{Definition}

\begin{Remark}
The semi-direct product is the most basic group extension of $Q$ by $K$ for any given outer action $\psi$ of $Q$ on $K$; $H^2(Q; Z(K))$ acts on the semi-direct product for the given outer action and permutes it to any other group extension. The semi-direct product for a given outer action $\psi$ is the only group extension that has an isomorphic copy of $Q$ living inside $G$ with an embedding $k: Q \rightarrow G$ satisfying $k \circ \sigma = id_Q$.
\end{Remark}

\begin{Theorem}[Normal Form for Semi-Direct Products] \label{normal-form}
Here we are viewing $K$ and $Q$ as subgroups of $G$. Let $G = K \rtimes Q$, where $K$ is generated by $\{\alpha_1, \alpha_2, \ldots\}$, $Q$ is generated by $\{\beta_1, \beta_2, \ldots\}$, and $\psi: Q \rightarrow Out(K)$ is the outer action. Then each element of $G$ admits a normal form as a product of generators of $Q$ and $K$ where all the generators of $K$ are on the left and all the generators of $Q$ are on the right.
\end{Theorem}

\begin{Proof}
Proof is by double induction on the number of generators of $Q$ not all on the right and the number of generators of $K$ between the collection of generators of $Q$ entirely on the right and the first generator of $Q$ not entirely on the right in a given representation of $g \in G$.

\indexspace

First, suppose there is one generator of $Q$ separated from the collection of generators of $Q$ entirely on the right and there is one generator of $K$ separating this generator of $Q$ from the collection of generators of $Q$ entirely on the right. Then we have $g = \alpha_{i_1} \cdot \ldots \cdot \alpha_{i_p} \cdot \beta_j \cdot \alpha_i \cdot \beta_{j_1} \cdot \ldots \cdot \beta_{j_q}$. Then $\beta_j \cdot \alpha_i \cdot \beta_j^{-1} = \psi(\beta_j)(\alpha_i)$, which implies $\beta_j \cdot \alpha_i = \psi(\alpha)(\beta_j) \cdot \beta_j$. Thus $g = \alpha_{i_1} \cdot \ldots \cdot \alpha_{i_p} \cdot \psi(\alpha)(\beta_j) \cdot \beta_j \cdot \beta_{j_1} \cdot \ldots \cdot \beta_{j_q}$, and $g$ admits a normal form. This establishes the basis case for the second induction.

\indexspace

Suppose the inductive hypothesis for the second induction; that is, suppose whenever there is one generator of $Q$ separated from the collection of generators of $Q$ entirely on the right and there are $m-1$ generators of $K$ separating this generator of $Q$ from the collection of generators of $Q$ entirely on the right, then this element of $G$ admits a normal form where all the generators of $K$ come first and all the generators of $Q$ come last. Now, suppose there is one generator of $Q$ separated from the collection of generators of $Q$ entirely on the right and there are $m$ generators of $K$ separating this generator of $Q$ from the collection of generators of $Q$ entirely on the right. Then we have $g = \alpha_{i_1} \cdot \ldots \cdot \alpha_{i_p} \cdot \beta_{j_1} \cdot \ldots \cdot \beta_{j_m} \cdot \alpha_i \cdot \beta_{j_{m+1}} \cdot \ldots \cdot \beta_{j_{m+q}}$. Then $\beta_{j_m} \cdot \alpha_i \cdot \beta_{J_m}^{-1} = \psi(\beta_{j_m})(\alpha_i)$, which implies $\beta_{j_m} \cdot \alpha_i = \psi(\beta_{j_m})(\alpha_i) \cdot \beta_{j_m}$. Thus $g = \alpha_{i_1} \cdot \ldots \cdot \alpha_{i_p} \cdot \beta_{j_1} \cdot \ldots \cdot \beta_{j_{m-1}} \cdot \psi(\alpha_i)(\beta_{j_m}) \cdot \beta_{j_m} \cdot \ldots \cdot \beta_{j_{m+q}}$, and $g$ admits a normal form. This establishes the inductive step for the second induction.

\indexspace

This also establishes the base case for the first induction.

\indexspace

Now, suppose the inductive hypothesis for the first induction; that is, suppose that whenever $g \in G$ has a representation where there are $n-1$ generators of $Q$ separated from the collection of generators of $Q$ entirely on the right and there are $m$ generator of $K$ separating the generator of $Q$ closest to the collection of generators of $Q$ entirely on the right, then this element of $G$ admits a representation with $n-1$ generators of $Q$ separated from the collection of generators of $Q$ entirely on the right and there are $(m-1) + l$ generators of $K$ separating the generator of $Q$ closest to the collection of generators of $Q$ entirely on the right for some natural number $l$. Suppose that $g \in G$ has a representation where there are $n$ generators of $Q$ separated from the collection of generators of $Q$ entirely on the left and there are $m$ generator of $K$ separating the generator of $Q$ closest to the collection of generators of $Q$ entirely on the left. Then we have $g = \alpha_{i_1} \cdot \ldots \cdot \alpha_{i_p} \cdot \beta_{j_1} \cdot \ldots \cdot \beta_{j_m} \cdot \alpha_i \cdot \beta_{j_{m+1}} \cdot \ldots \cdot \beta_{j_{m+q}}$. Then $\beta_{j_m} \cdot \alpha_i \cdot \beta_{j_m}^{-1} = \psi(\beta_{j_m})(\alpha_i)$, which implies $\beta_{j_m} \cdot \alpha_i = \psi(\beta_{j_m})(\alpha_i) \cdot \beta_{j_m}$. Thus $g = \alpha_{i_1} \cdot \ldots \cdot \alpha_{i_p} \cdot \beta_{j_1} \cdot \ldots \cdot \beta_{j_{m-1}} \cdot \psi(\beta_{j_m})(\alpha_i) \cdot \beta_{j_m} \cdot \beta_{j_{m+1}} \cdot \ldots \cdot \beta_{j_{m+q}}$, and $g$ admits a normal form. 
\end{Proof}

\begin{Corollary}[Presentations for Semi-Direct Products]\label{corpres-semi-dir-products}
Let $G = K \rtimes Q$, where $K$ is presented by $\{\alpha_1, \alpha_2, \ldots, | r_1, r_2, \ldots \}$ and $Q$ is presented by $\{\beta_1, \beta_2 \ldots, | s_1, s_2, \ldots \}$. Then $G$ admits a presentation as

\begin{equation}
G \cong \langle \alpha_i, \beta_j | r_k, s_l, \beta_j\alpha_i(\psi(\beta_j)(\alpha_i))^{-1} \rangle 
\end{equation}

where $\psi$ is the outer action and each $\psi(\beta_i)(\alpha_j)$ is a word in the $\alpha_j$'s 
\end{Corollary}

\begin{Proof}
Since by Corollary \ref{normal-form}, each word has a unique normal form given by sliding all the $\alpha$'s to the left and all the $\beta$'s to the right, the slide relators, as well as the defining relators from each group, are relators in the semi-direct product. Since each word can be put in normal form using only these relators, there is a presentation using only these relators.
\end{Proof}

\begin{Remark}
A relator $\beta_i\alpha_j\beta_i^{-1}(\psi(\beta_i)(\alpha_j))^{-1}$ is sometimes called a \textbf{slide relator}, and the word $\psi(\beta_i)(\alpha_j)$ represents ``the price for moving $\beta_i$ across $\alpha_j$''. Note that in a direct product, we have the trivial outer action of $Q$ on $K$, and so we have $\psi(\beta_i)(\alpha_j) = \alpha_j$; there is, in some sense, no price to pay for sliding an element of the quotient group across an element of the kernel group - the two groups commute. Semi-direct products are like direct products, except that there is a price to pay for sliding a quotient group element across a kernel group element.
\end{Remark}
        
    %
    %

    \newchapter{HANDLEBODY REVERSE TO PLUS}{A HANDLEBODY-THEORETIC REVERSE TO THE PLUS CONSTRUCTION}{A GEOMETRIC REVERSE TO THE PLUS CONSTRUCTION IN HIGH DIMENSIONS}
    \label{sec:LabelForChapter3}

     
        
    
        \section[A Handlebody-Theoretic Reverse to the Plus Construction]{A Handlebody-Theoretic Reverse to the Plus Construction}
        \label{sec:LabelForChapter3:Section1}
        In this section, we will describe our partial solution to the Reverse Plus Problem. Our solution applies to superperfect (defined in Definition \ref{defsuperperfect} below), finitely presented kernel groups. Also, our solution applies to the case that the total group $G$ of the group extension $1 \rightarrow K \rightarrow G \rightarrow Q \rightarrow 1$ is a semi-direct product (defined in Definition \ref{defsemi-dir-prod} above). This is an important special case of a hard problem.

\indexspace

It is, however, we believe, easy to use and easy to understand. We feel that in the situations where our solution applies (superperfect, finitely presented kernel group and semi-direct product for the total group of the group extension), we have reduced the topological problem of solving the Reverse Plus Problem to an algebraic problem of computing a semi-direct product of two groups $Q$ and $K$ by identifying an outer action of $Q$ on $K$; this is supposed to be the goal of algebraic topology in general.

\indexspace

\begin{Definition} \label{defsuperperfect}
A group $G$ is said to be \textbf{superperfect} if its first two homology groups are 0, that is, if $H_1(G) = H_2(G) = 0$. (Recall a group is \textbf{perfect} if and only if its first homology group is 0.)
\end{Definition}

\begin{Example}
A perfect group is superperfect if it admits a finite, \textit{balanced} presentation, that is, a finite  presentation with the same number of generators as relators. (The converse for finitely presented superperfect groups is false.)
\end{Example}

\begin{Lemma} \label{lemsphere-elts} 
Let $S$ be a superperfect group. Let $K$ be a cell complex which has fundamental group isomorphic to $S$. Then all elements of $H_2(K)$ can be killed by attaching 3-cells.
\end{Lemma}

\begin{Proof}
By Proposition 7.1.5 in \cite{Geoghegan}, there is a $K(S,1)$ which is formed from $K$ by attaching cells of dimension 3 and higher. Let $L$ be such a $K(S,1)$. Then $L^3$ is formed from $K^2$ by attaching only 3-cells, and $H_2(L^3) \cong H_2(L)$, as $L$ is formed from $L^3$ by attaching cells of dimension 4 and higher, which cannot affect $H_2$. But $H_2(L) \cong H_2(S)$ by definition and $H_2(S) \cong 0$ by hypothesis. Thus, all elements of $H_2(K)$ can be killed by attaching 3-cells.
\end{Proof}

\begin{Lemma}[Equivariant Attaching of Handles] \label{lemequiv-attach}
Let $M^n$ be a smooth manifold, $n \ge 5$, with $M$ one boundary component of $W$ with $\pi_1(M) \cong G$. Let $P \unlhd G$ and $Q = G/P$. Let $\overline{M}$ be the cover of $M$ with fundamental group P and give $H_*(\overline{M}; \mathbb{Z})$ the structure of a $\mathbb{Z}Q$-module. Let $2k + 1 \le n$ and let $S$ be a finite collection of elements of $H_k(M; \mathbb{Z})$  which all admit embedded spherical representatives which have trivial tubular neighborhoods. If $k = 1$, assume all elements of $S$ represent elements of $P$.

\indexspace

Then one can \textbf{equivariantly attach $(k+1)$-handles across S}, that is, if $\overline{S} = \{ s_{j,q}\ |\ q \in Q\}$ is the collection of lifts of elements of $S$ to $\overline{M}$, one can attach $(k+1)$-handles across tubular neighborhoods of the $s_{j,q}$ so that each lift $s_{j,q}$ projects down via the covering map $p$ to an element $s_j$ of $S$ and so that the covering map extends to send each $(k+1)$-handle $H_{j,q}$ attached across a tubular neighborhood of $s_{j,q}$ in $\overline{M}$ bijectively onto a handle attached across the projection via the covering map of the tubular neighborhood of the element $s_j$ in $M$.
\end{Lemma}

\begin{Proof}
$H_k(\overline{M}; \mathbb{Z})$ has the structure of a $\mathbb{Z}Q$-module. The action of $Q$ on $\overline{S}$ permutes the elements of $S$. For each embedded sphere $s_j$ in $S$, lift it via its inverse images under the covering map to a pairwise disjoint collection of embedded spheres $s_{j,q}$. (This is possible since a point of intersection or self-intersection would have to project down to a point of intersection or self-intersection (respectively) by the evenly-covered neighborhood property of covering spaces.) The $s_{j,q}$ all have trivial tubular neighborhoods. Attach an $(k+1)$-handle across the tubular neighborhood of the elements $s_j$ of the $S$. For each $j \in \{1, \ldots, |S|\}$ and $q \in Q$ attach an $(k+1)$-handle across the spherical representative $s_{j,q}$; extend the covering projection so it projects down in a bijective fashion from the handle attached along $s_{j,q}$ onto the handle we attached along $s_j$.
\end{Proof}

\begin{Lemma} \label{lemkernel}
Let $A, B,$ and $C$ be $R$-modules, with $B$ a free $R$-module (on the basis $S$), and let $\Theta: A \bigoplus B \rightarrow C$ be an $R$-module homomorphism. Suppose $\Theta|_A$ is onto. Then $\ker(\Theta) \cong \ker(\Theta|_A) \bigoplus B$.
\end{Lemma}

\begin{Proof}
Define $\phi: \ker(\Theta|_A) \bigoplus B \rightarrow \ker(\Theta)$ as follows. For each $s \in S$, where $S$ is a basis for $B$, choose $\alpha(s) \in A$ with $\Theta(\alpha(s),0) = \Theta(0,s)$, as $\Theta|_A$ is onto. Extend $\alpha$ to a homomorphism from $B$ to $A$, and note that $\alpha$ has the same property for all $b \in B$. Then set $\phi(x,b) = (x-\alpha(b),b)$.

\indexspace

(Well-defined) Let $x \in \ker(\Theta|_A)$ and $b \in B$. Then $\Theta(\phi(x,b)) = \Theta(x-\alpha(b),b) = \Theta(x,0) + \Theta(-\alpha(b),0) + \Theta(0,b) = 0 + -\Theta(\alpha(b),0) + \Theta(0,b) = 0 + -\Theta(0,b) + \Theta(0,b) = 0$. So, $\phi$ is well-defined.

\indexspace	

Define $\psi: \ker(\Theta) \rightarrow \ker(\Theta|_A) \bigoplus B$ by $\psi(z) = (\pi_1(z)+\alpha(\pi_2(z)),\pi_2(z))$, where $\pi_1: A \bigoplus B \rightarrow A$ and $\pi_2: A \bigoplus B \rightarrow B$ are the canonical projections.

\indexspace

(Well-defined) Let $z \in \ker(\Theta)$. It is clear that $\pi_2(z) \in B$, so it remains to prove that $\pi_1(z) + \alpha(\pi_2(z)) \in \ker(\Theta|_A)$. [Note $\Theta(z) = \Theta|_A(\pi_1(z)) + \Theta|_B(\pi_2(z)) \Rightarrow \Theta|_A(\pi_1(z)) = -\Theta|_B(\pi_2(z))$. Note also, by definition of $\alpha$, $\Theta(\alpha(\pi_2(z))) = \Theta(0,\pi_2(z))$]. We compute $\Theta|_A(\pi_1(z) + \alpha(\pi_2(z))) = \Theta|_A(\pi_1(z)) + \Theta(\alpha(\pi_2(z)))$ = $-\Theta|_B(\pi_2(z)) + \Theta(0,\pi_2(z)) = -\Theta(0,\pi_2(z)) + \Theta(0,\pi_2(z)) = 0$. So, $\psi$ is well-defined.

\indexspace

(Homomorphism) Clear.

\indexspace

(Inverses) Let $(x,b) \in \ker(\Theta|_A) \bigoplus B$. The $\psi(\phi(x,b)) = \psi(x-\alpha(b),b) = (\pi_1(x-\alpha(b),b)+\alpha(\pi_2(x-\alpha(b),b)),\pi_2(x-\alpha(b),b)) = (x-\alpha(b)+\alpha(b),b) = (x,b)$.

\indexspace

Let $z \in \ker(\Theta)$. Then $\phi(\psi(z)) = \phi(\pi_1(z)+\alpha(\pi_2(z)),\pi_2(z)) = (\pi_1(z)+\alpha(\pi_2(z))-\alpha(\pi_2(z)),\pi_2(z)) = (\pi_1(z),\pi_2(z)) = z$.

\indexspace

So, $\phi$ and $\psi$ are inverses of each other, and the lemma is proven.
\end{Proof}

\begin{Definition}
A $k$-handle is said to be \textbf{trivially attached} if and only if it is possible to attach a canceling $k+1$-handle.
\end{Definition}

Here is our solution to the Reverse Plus Problem in the high-dimensional manifold category.

\begin{RestateTheorem}{Theorem}{thmsemi-s-cob}{An Existence Theorem for Semi-s-Cobordisms}
Given $1 \rightarrow S \rightarrow G \rightarrow Q \rightarrow 1$ where $S$ is a finitely presented superperfect group, $G$ is a semi-direct product of $Q$ by $S$, and any closed $n$-manifold $N$ with $n \ge 6$ and $\pi_1(N) \cong Q$, there exists a solution $(W, N, N_-)$ to the Reverse Plus Problem for which $N \hookrightarrow W$ is a simple homotopy equivalence.
\end{RestateTheorem}

\begin{Proof}
Start by taking $N$ and crossing it with $\mathbb{I}$. Let $Q \cong  \langle \alpha_1, \ldots, \alpha_{k_1} | r_1, \ldots, r_{l_1} \rangle$ be a presentation for $Q$. Let $S \cong  \langle \beta_1, \ldots, \beta_{k_2} | s_1, \ldots, s_{l_2} \rangle $ be a presentation for $S$. Take a small $n$-disk $D$ inside of $N \times \{1\}$. Attach a trivial 1-handle $h^1_i$ for each $\beta_i$ in this disk $D$. Note that because they are trivially attached, there are canceling 2-handles $k_i^2$, which may also be attached inside the disk together with the 1-handles $D \cup \{h^1_i\}$. We identify these 2-handles now, but do not attach them yet. They will be used later.

\indexspace

Attach a 2-handle $h^2_j$ across each of the relators $s_j$ of the presentation for $S$ in the disk together with the 1-handles $D \cup \{h^1_i\}$, choosing the framing so that it is trivially attached in the manifold that results from attaching $h_i^1$ and $k_i^2$ (although we have not yet attached the handles $k_i^2$). Note that because they are trivially attached, there are canceling 3-handles $k_j^3$, which may also be attached in the portion of the manifold consisting of the disk $D$ together with the 1-handles $\{h^1_i\}$ and the 2-handles $\{k^2_i\}$. We identify these 3-handles now, but do not attach them yet. They will be used later.

\indexspace

By Corollary \ref{corpres-semi-dir-products}, let 
\begin{equation}
\begin{split}
G \cong & \langle \alpha_1, \ldots, \alpha_{k_1}, \beta_1, \ldots, \beta_{k_2} | r_1, \ldots, r_{l_1}, s_1, \ldots, s_{l_2}, \beta_j\alpha_i(\psi(\beta_j)(\alpha_i))^{-1} \rangle 
\end{split}
\end{equation}

be a presentation for $G$. Attach a 2-handle $f^2_{i,j}$ for each relator $q_jk_iq_j^{-1}\phi(k_i)^{-1}$, choosing the framing so that it is trivially attached in the result of attaching the $h_i^1$, $k_i^2$, $h_2^j$ and $k_j^3$. This is possible since each of the relators becomes trivial when the $k_i^2$'s and $k_i^3$'s are attached. Note that because the $f^2_{i,j}$ are trivially attached, there are canceling 3-handles $g_{i,j}^3$. We identify these 3-handles now, but do not attach them yet. They will be used later. Call the resulting cobordism with only the $h_i^1$'s, the $h_j^2$'s, and the $f_{i,j}^2$'s attached $(W', N, M')$ and call the right-hand boundary $M'$.

\indexspace

Note that we now have $\pi_1(N) \cong Q$, $\pi_1(W') \cong G$, and $\iota_\#: \pi_1(M') \rightarrow \pi_1(W)$ an isomorphism because, by inverting the handlebody decomposition, we are starting with $M'$ and adding $(n-1)$- and $(n-2)$-handles, which do not affect $\pi_1$ as $n \ge 6$.

\indexspace

Consider the cover $\overline{W'}$ of $W'$ corresponding to $S$. Then the right-hand boundary of this cover, $\overline{M'}$, also has fundamental group isomorphic to $S$ by covering space theory. Also, the left-hand boundary of this cover, $\widetilde{N}$, has trivial fundamental group.

\indexspace

Consider the handlebody chain complex $C_*(\overline{W'}, \widetilde{N}; \mathbb{Z})$. This is naturally a $\mathbb{Z}Q$-module complex. It looks like

\begin{diagram}[size=14.5pt]
0 & \rTo & C_3(\overline{W'}, \widetilde{N}; \mathbb{Z}) & \rTo & C_2(\overline{W'}, \widetilde{N}; \mathbb{Z}) & \rTo^{\partial} & C_1(\overline{W'}, \widetilde{N}; \mathbb{Z}) & \rTo & C_0(\overline{W'}, \widetilde{N}; \mathbb{Z}) & \rTo & 0 \\
= &      & =                                            &       & \cong                                 
                                              &                 & \cong
                                              &      & =      
                           &      & = \\
0 & \rTo & 0                     & \rTo & \bigoplus_{i=1}^{l_2} \mathbb{Z}Q 
\oplus \bigoplus_{j=1}^{k_1*	k_2} \mathbb{Z}Q & \rTo^{\partial} & 
\bigoplus_{i=1}^{k_2} \mathbb{Z}Q & \rTo & 0                     & \rTo & 0 \\
\end{diagram}

\indexspace

where $C_2(\overline{W'}, \widetilde{N}; \mathbb{Z})$ decomposes as $A = \bigoplus_{i=1}^{l_2} \mathbb{Z}Q$, which has a $\BZ Q$-basis obtained by arbitrarily choosing one lift of the 2-handles for each of the $h^2_j$, and $B = \bigoplus_{j=1}^{k_1 \cdot k_2} \mathbb{Z}Q$, which has a $\BZ Q$-basis obtained by arbitrarily choosing one lift of the 2-handles for each of the $f^2_{i,j}$. Set $C = C_1(\overline{W'}, \widetilde{N}; \mathbb{Z}) \cong \bigoplus_{i=1}^{k_2} \mathbb{Z}Q$ (as $\BZ Q$-modules). Choose a preferred basepoint $\overline{*}$ and a preferred lift of the the disk $D$ to a disk $\overline{D}$ in $\overline{M}$. Decompose $\partial$ as $\partial_{2,1} = \partial|_A$ and $\partial_{2,2} = \partial|_B$

\indexspace

Since $S$ is perfect, we must have $l_2 \ge k_2$, as we must have as many or more relators as we have generators in the presentation for $S$ to have no 1-dimensional homology. 

\indexspace

We examine the contribution of $\partial_{2,1}$ to $H_2(\overline{W'}, \widetilde{N}; \mathbb{Z})$. It will be useful to first look downstairs at the $\BZ$-chain complex for $(W', N)$. Let $A'$ be the submodule of $C_2(W',N;\mathbb{Z})$ determined by the $h^2_j$'s and let $C'$ be $C_1(W',N;\BZ)$, which is generated by the $h^1_i$'s. Then $A'$ is a finitely generated free abelian group, so, the kernel $K'$ of $\partial'_{2,1}: A' \rightarrow C'$ is a subgroup of a finitely generated free abelian group, and thus $K'$ is a finitely generated free abelian group, say on the basis $\{k_1, \ldots, k_a\}$

\indexspace

\begin{Claim}
$\ker(\partial_{2,1})$ is a free $\mathbb{Z}Q$-module on a generating set of cardinality $|a|$.
\end{Claim}

\begin{Proof}
The disk $D$ has $|Q|$ lifts of itself to $\overline{M}$, ala Lemma \ref{lemequiv-attach}. Now, $Q$ acts as deck transformations on $\overline{M}$, transitively permuting the lifts of $D$ as the cover $\overline{M}$ is a regular cover. A preferred basepoint $\overline{*}$ and a preferred lift of the the disk $D$ to a disk $\overline{D}$ in $\overline{M}$ have already been chosen for the identification of $C_*(\overline{W'}, \widetilde{N}; \BZ)$ with the $\BZ Q$-module $C_*(W', N; \mathbb{Z}Q)$. Let the handles attached inside the preferred lift $\overline{D}$ be our preferred lifts $\overline{h^1_i}$ and let the lifts of the $h_j^2$s that attach to $\overline{D} \cup (\cup \overline{h}_i^1)$ be our preferred 2-handles $\overline{h}_j^2$.

\indexspace

Note that none of the $q\overline{h^1_i}$ spill outside the disk $q\overline{D}$ and none of the $q\overline{h^2_j}$ spill outside the disk $q\overline{D} \cup (\cup q\overline{\overline{h}^1_i})$. This implies $\partial_{2,1}(\overline{\overline{h}^2_j}) \in \{z_ih_i^1\ |\ z_i \in \BZ\} \le \{z_i q_i h_i^1 \ |\ z_iq_i \in \BZ Q \}$ and so $\partial_{2,1}(q\overline{h^2_j}) \in \{z_iqh_i^1\ |\ z_i \in \BZ q_t \le \BZ Q \}$. This mean if $q_1 \neq q_2$ are in $Q$ and $\overline{c_1}$ and $\overline{c_2}$ are lifts of chains in $A$ to $\overline{D}$, then $\partial_{2,1}(q_1\overline{c_1}  + q_2\overline{c_2}) = 0 \in \BZ Q$ if and only if $\partial_{2,1}(\overline{c_1}) = \partial_{2,1}(\overline{c_2}) = 0 \in \BZ$; ($\ddag$).

\indexspace

With this in mind, let $\overline{k_i}$ be a lift of the chain $k_i$ in a generating set for $K'$ in the disk $D$ to $\overline{D}$. Then $\partial_{2,1}(\overline{k_i}) = 0$. Moreover, $Q$ transitively permutes each $\overline{k_i}$ with the other lifts of $k_i$ to the other lifts of $D$. Now, suppose $\partial_{2,1}(c) = 0$, with $c$ an element of $C_2(W', N; \BZ Q)$. By ($\ddag$), we must have $c = \sum_{t=1}^m n_tq_t\overline{k_t}$ with $n_t \in \BZ$ and $q_t \in Q$. This proves the $\overline{k_t}$'s generate $\ker(\partial_{2,1})$.

\indexspace

Finally, suppose some linear combination $\sum_{i=1}^a (\sum n_tq_t)\overline{k_i}$ is zero. Then, as  $q_{t_1}\overline{k_{t_1}}$ and $q_{t_2}\overline{k_{t_2}}$ cannot cancel if $q_{t_1} \neq q_{t_2}$, it follows that all $n_t$ are zero. This proves the $\overline{k_i}$'s are a free $\BZ Q$-basis for $\ker(\partial_{2,1})$. This proves the claim.
\end{Proof}

\indexspace

Now, we have $\partial_2: A \bigoplus B \rightarrow C$. Recall $S$ is a finitely presented, superperfect group, and W' contains a 1-handle for each generator and a 2-handle for each relator in a chosen finite presentation for $S$. It then follows that $\ker(\partial_1)/\text{im}(\partial_2|_A) \cong 0$, as if $\Lambda$ contains the collection of lifts of 1-handles for each generator of $S$ and the collection of lifts of 2-handles for each relator of $S$, then $\Lambda = 0$ as a $\BZ Q$-modules and $\Lambda = \ker(\partial_1)/\text{im}(\partial_2|_A)$. But $C_0(\overline{W'}, \widetilde{N}; \mathbb{Z}) = 0$, so  $\ker(\partial_1) = C$. This implies $\partial_2|_A$ is onto. By Lemma \ref{lemkernel}, we have that $\ker(\partial_2) \cong  \ker(\partial_2|_A) \bigoplus B$. By the previous claim, $\ker(\partial_2|_A)$ is a free and finitely generated $\mathbb{Z}Q$-module. Clearly, $B$ is a free and finitely generated $\mathbb{Z}Q$-module. Thus, $\ker(\partial_2) \cong H_2(\overline{W'}, \widetilde{N}; \mathbb{Z})$ is a free and finitely generated $\mathbb{Z}Q$-module.

\indexspace

By Lemma \ref{lemsphere-elts}, we may choose spherical representatives for all elements of $H_2(\overline{W'}; \mathbb{Z})$. By the Long Exact Sequence in homology for $(\overline{W'}, \widetilde{N})$, we have 

\begin{diagram}
\cdots & \rTo & H_2(\overline{W'};\mathbb{Z}) & \rTo   & H_2(\overline{W'}, \widetilde{N}; \mathbb{Z}) & \rTo  & H_1(\widetilde{N}; \mathbb{Z}) & \rTo & \cdots \\
       &      & =                             &        & =                          
                           &       & \cong                          &      & \\
\cdots & \rTo & H_2(\overline{W'};\mathbb{Z}) & \rOnto & H_2(\overline{W'}, \widetilde{N}; \mathbb{Z}) & \rTo & 0                               & \rTo & \cdots \\
\end{diagram}

so any element of $H_2(\widetilde{W'}, \overline{N}; \mathbb{Z})$ also admits a spherical representative. 

\indexspace

So, we may choose spherical respresentatives for any element of $H_2(W', N; \BZ Q)$. Let $\{s_k\}$ be a collection of embedded, pair-wise disjoint 2-spheres which form a free, finite $\BZ Q$-basis for $H_2(W', N; \BZ Q)$.

\indexspace

Note that the $\{s_k\}$ can be arranged to live in right-hand boundary $M'$ of $W'$. To do this, view $W'$ upside-down, so that it has $(n-2)-$ and $(n-1)-$handles attached. For each $s_k$, make it transverse to the (2-dimensional) co-core of each $(n-2)-$handle, then blow it off the handle by using the product structure of the handle less the co-core; do the same thing with the $(n-1)-$handles. Finally, use the product structure of $N \times \BI$ to push $s_k$ into the right-hand boundary.

If we add the $k_i^2$, $h_j^3$ and $g_{i,j}^3$ to $W'$, and similarly make sure the $k_i^2$s, $k_j^3$s, and $g_{i,j}^3$s do not intersect the $\{s_k\}$s, and call the resulting cobordism $W''$, we can think of the $\{s_k\}$ as living in the right-hand boundary of $(W'', N, M'')$. Note that $W''$ is diffeomorphic to $N \times \BI$.

\indexspace

We wish to attach 3-handles along the collection $\{s_k\}$ and, later, 4-handles complimetary to those 3-handles. \textit{A priori}, this may be impossible; for instance, there is a framing issue. To make this possible, we borrow a trick from \cite{G-T3} to alter the 2-spheres to a useable collection without changing the elements of $H_2(W',N; \BZ Q)$ they represent.

\indexspace

\begin{Claim}
For each $s_k$, we may choose a second embedded 2-sphere $t_k$ with the property that 
\begin{itemize}
\item $t_k$ represents the same element of $\pi_2(M'')$ as $s_k$ (as elements of $\pi_2(W')$, they will be different)
\item each $t_k$ misses the attaching regions of all the $\{h_i^1\}, \{k_i^2\}, \{h_j^2\}, \{k_j^3\}, \{f_{i,j}^2\}$ and $\{g_{i,j}^3\}$
\item the collection of $\{t_k\}$ are pair-wise disjoint and disjoint from the entire collection $\{s_k\}$
\end{itemize}
\end{Claim}

\begin{Proof}
Note that each canceling (2,3)-handle pairs $h_j^2$ and $k_j^3$ and $f_{i,j}^2$ and $g_{i,j}^3$ form an $(n+1)$-disks attached along an $n$-disk which is a regular neighborhood of a 2-disk filling the attaching sphere of the 2-handle. Also, each canceling (1,2)-handle $h_i^1$ and $k_i^2$ forms an $(n+1)$-disk in $N \times \{1\}$ attached along an $n$-disk which is a regular neighborhood of a 1-disk filling the attaching sphere of the 1-handle. We may push a given $s_k$ off the (2,3)-handle pairs and then off the (1,2)-handle pairs, making sure not to pass back into the (2,3)-handle pairs. Let $t_k$ be the end result of the pushes. Make the collection $\{t_k\}$ pair-wise disjoint and disjoint from the $\{s_k\}$'s by tranversality, making sure not to pass back into the (1,2)- or (2,3)-handle pairs.
\end{Proof}

Replace each $s_k$ with $s_k\#(-t_k)$, an embedded connected sum of $s_k$ with a copy of $t_k$ with its orientation reversed.

\indexspace

Since the $t_k$'s miss \emph{all} the handles attached to the original collar $N \times \BI$, they can be pushed into the right-hand copy of $N$. Thus, $s_k$ and $s_k\#(-t_k)$ represent the same element of $H_2(W', N; \BZ Q)$. Hence, the collection $\{s_k\#(-t_k)\}$ is still a free basis for $H_2(W', N; \BZ Q)$. Furthermore, each $s_k\#(-t_k)$ bounds an embedded 3-disk in the boundary of $W''$. This means each $s_k\#(-t_k)$ has a product neighborhood structure, and we may use it as the attaching region for a 3-handle $h_l^3$. Choose the framing of $h_l^3$ so that it is a trivially attached 3-handle with respect to $W''$, and choose a canceling 4-handle $k_l^4$. We identify these 4-handles now, but do not attach them yet. They will be used later. Call the resulting cobordism with the $h_i^1$, $h_j^2$, $f_{i,j}^2$, and $h_l^3$ attached $(W''', N, M)$. Let $W^{(iv)}$ be $M \times \mathbb{I}$ with the $k^2_i$, $k^3_j$, and $k^4_k$'s attached. Then $W''' \bigcup_{M} W^{(iv)}$ has all canceling handles and so is diffeomorphic to $N \times \mathbb{I}$. Clearly, $W''' \bigcup_M W^{(iv)}$ strong deformation retracts onto the right-hand boundary $N$. Despite all the effort put into creating $(W''', M, N)$, $(W^{(iv}, M, N)$, or, more precisely, $(W^{(iv}, N, M)$ (modulo torsion) will be seen to satisfy the conclusion of the theorem.

\indexspace

We are note yet finished with $(W''', N, M)$ yet. In order to prove $(W^{(iv}, M, N)$ satisfies the desired properties, we must study $W'''$ more carefully. Note that since $\ker(\partial_2)$ is a free, finitely generated $\BZ Q$-module and $\{h_k^3\}$ is a set whose attaching spheres are a free $\BZ Q$-basis for $\ker(\partial_2)$, $\partial_3: C_3(W''', N; \BZ Q) \rightarrow C_2(W''', N; \BZ Q)$ is onto and has no kernel. This means $H_3(W''', N; \BZ Q) \cong 0$. Clearly, $H_*(W''', N; \BZ Q) \cong 0$ for $* \ge 4$ as $C_*(W''', N; \BZ Q) \cong 0$ for $* \ge 4$.

\indexspace

Thus, $H_*(\overline{W'''}, \widetilde{N}; \mathbb{Z}) \cong 0$, i.e., $H_*(W''', N; \BZ Q) \cong 0$. (*)

\indexspace

However, this is not the only homology complex we wish to prove acyclic; we also wish to show that $H_*(W''', M; \BZ Q) \cong 0$.By noncompact Poincare duality, we can do this by showing that the relative cohomology with compact supports is 0, i.e., $H^*_c(\overline{W'''}, \widetilde{N}; \BZ) \cong 0$.

%
%
%
%
%
%
%
%
%

\indexspace

By the cohomology with compact supports, we mean to take the chain complex that has linear functions $f: C_i(\overline{W'''}, \widetilde{N}; \BZ) \rightarrow \BZ$ from the relative handlebody complex of the intermediate cover of $W'''$ with respect to $K$ to $\BZ$ relative to $\widetilde{N}$, that is, that sends all of the handles of the universal cover of $N$ to 0 and that is nonzero on only finitely many of the $qh_j$'s. The fact that $\delta$ is not well-defined, that is, that $g$ has compact supports depends on the fact that $C_*(\overline{W'''}, \widetilde{N}; \BZ)$ is locally finite, which in turn depends on the fact that $\overline{W'''}$ is a covering space of a compact manifold, with finitely many handles attached.The co-boundary map $\delta_*$ will send a cochain $f$ in $C^i_c(\overline{W'''}, \widetilde{N}; \BZ)$ to the cochain $g$ in $C^{i+1}_c(\overline{W'''}, \widetilde{N}; \BZ)$ which sends $g(\partial(n_jq_ih_j)$ to $\delta(f)(n_jq_ih_j)$ for $q_i \in Q$ and $n_jh_j \in C_i(W''', N; \BZ)$. 

\indexspace

%
%
%

Clearly, $\delta_1: C^0_c(\overline{W'''}, \widetilde{N}; \BZ) \rightarrow C^1_c(\overline{W'''}, \widetilde{N}; \BZ)$ and $\delta_4: C^3_c(\overline{W'''}, \widetilde{N}; \BZ) \rightarrow C^4_c(\overline{W'''}, \widetilde{N}; \BZ)$ are the zero maps. This means we must show $ker(\delta_2) = 0$, i.e., $\delta_2$ is 1-1, and $im(\delta_3) = C_3$, i.e., $\delta_3$ is onto. Finally, we must show exactness at $C^2_c$, that is, we must show $im(\delta_2) = ker(\delta_3)$.

Consider the acyclic complex

\begin{diagram}[size=14.5pt]
0 & \rTo & C_3(\overline{W'''}, \widetilde{N}; \mathbb{Z}) & \rTo & C_2(\overline{W'''}, \widetilde{N}; \mathbb{Z}) & \rTo^{\partial} & C_1(\overline{W'''}, \widetilde{N}; \mathbb{Z}) & \rTo & 0 \\
\end{diagram}

This admits a section $\iota: C_1(\overline{W'''}, \widetilde{N}; \mathbb{Z}) \rightarrow C_2(\overline{W'''}, \widetilde{N}; \mathbb{Z})$ with the property that $\partial_3(C_3) \bigoplus \iota(C_1) = C_2$

\indexspace

($ker(\delta_2) = 0$) Let $f \in C^1_c(\overline{W'''}, \widetilde{N}; \BZ)$ be non-zero, that is, let $f: C_1(\overline{W'''}, \widetilde{N}; \BZ) \rightarrow 0$ have compact support and that there is a $c_1 \in C_1(\overline{W'''}, \widetilde{N}; \BZ)$ with $c_1 \ne 0$ and $f(c_1) \ne 0$. As $\partial_2$ is onto, choose $c_2 \in C_2(\overline{W'''}, \widetilde{N}; \BZ)$ with $c_2 \ne 0$ and $\partial_2(c_2) = c_1$. The $\delta_2(f)(c_2) = f(\partial_2(c_2)) = f(c_1) \ne 0$, and $\delta_2(f)$ is not the zero cochain.

\indexspace

($im(\delta_3) = C^3$) Let $g \in C^3_c(\overline{W'''}, \widetilde{N}; \BZ)$ be a basis element with $g(qh^3_i) = 1$ and all other $g(q'h^3_{i'}) = 0$. We must show there is an $f \in C^2_C(\overline{W'''}, \widetilde{N}; \BZ)$ with $\delta_3(f) = g$. Consider $\partial_3(qh^3_i)$. This is a basis element for $C_2(\overline{W'''}, \widetilde{N}; \mathbb{Z})$.

Choose $f_{k,l} \in C^2_c(\overline{W'''}, \widetilde{N}; \BZ)$ to have $f_{k,l}(\partial_3(qh^3_i)) = 1$ and 0 otherwise. Then $\delta_3(f)(q_ih^3_j) = f(\partial_3(q_ih^3_j)) = 1 = g(qh^3_i)$.

\indexspace

This proves $\delta_3(f) = g$, and $\delta_3$ is onto.

\indexspace

($im(\delta_2) = ker(\delta_3)$) 

\indexspace

Clearly, if $f \in im(\delta_2)$, then $\delta_3(f) = 0$, as $\delta$ is a chain map.

\indexspace

Suppose $\delta_3(f) = 0$ but $f \ne 0$. Consider $\iota(qh^1_i) = c_{2,i} \in C_2(\overline{W'''}, \widetilde{N}; \mathbb{Z})$. This is a basis element for $C_2(\overline{W'''}, \widetilde{N}; \mathbb{Z})$.

\indexspace

Set $g(qh^1_i) = f(c_{2,i})$.

\indexspace

Then $\delta_2(g)(c_{2,i}) = g(\partial_2(c_{2,j})) = g(qh^1_i)  = f(c_{2,i})$, and we are done.

\indexspace

So, $H^*_C(\overline{W'''}, \widetilde{N}; \BZ) \cong 0$, so $H_*(\overline{W'''}, \overline{M}; \BZ) \cong 0$ by Theorem 3.35 in \cite{Hatcher}, and $H_*(W''', M; \BZ) \cong 0$

\indexspace

Note that we again have $\pi_1(N) \cong Q$, $\pi_1(W''') \cong G$, and $\iota_\#: \pi_1(M) \cong \pi_1(W''')$ an isomorphism, as attaching 3-handles does not affect $\pi_1$, and, dually, attaching $(n-3)$-handles does not affect $\pi_1$ for $n \ge 6$.

\indexspace

We read $W^{(iv)}$ from right to left. This is (almost) the cobordism we desire. (We will need to deal with torsion issues below.) Note that the left-hand boundary of $W^{(iv)}$ read right to left is $N$ and the right-hand boundary of $W^{(iv)}$ read right to left is $M$. Moreover, $W^{(iv)}$ read right to left is $N \times \BI$ with $[(n+1)-4]$-, $[(n+1)-3]$-, and $[(n+1)-2]$-handles attached to the right-hand boundary. Since $n \ge 6$, adding these handles does not affect $\pi_1(W^{(v)})$. Thus, we have $\iota_\#: \pi_1(N) \rightarrow \pi_1(W^{(v)})$ is an isomorphism; as was previously noted, $\pi_1(M) \cong G$.

\indexspace

Let $H: W''' \bigcup_M W^{(iv)} \rightarrow W''' \bigcup_M W^{(iv)}$ a strong deformation retraction onto the right-hand boundary $N$. We will produce a retraction $r: W''' \bigcup_M W^{(iv)} \rightarrow W^{(iv)}$. Then $r \circ H$ will restrict to a strong deformation retraction of $W^{(iv)}$ onto its right-hand boundary $N$. This, in turn, will yield a strong deformation retraction of $W^{(iv)}$ read right to left onto its left-hand boundary $N$.

\indexspace

Note that by (*), $H^*_C(\overline{W'''}, \widetilde{N}; \mathbb{Z}) \cong 0$. By Theorem 3.35 in \cite{Hatcher}, we have that $H_*(\overline{W'''}, \overline{M}; \mathbb{Z}) \cong 0$, and $H_*(W''', M; \mathbb{Z}Q) \cong 0$, respectively, by the natural $\BZ Q$ structure on $C_*(\overline{W'''}; \BZ)$.

\indexspace

To get the retraction $r$, we will use the following Proposition from \cite{G-T4}.

\begin{Proposition} \label{G-T4-lem42}
Let $(X,A)$ be a CW pair for which $A \hookrightarrow X$ induces a $\pi_1$ isomorphism. Suppose also that $L \unlhd \pi_1(A)$ and $A \hookrightarrow X$ induces $\BZ [\pi_1(A)/L]$-homology isomorphisms in all dimensions. Next suppose $\alpha_1, \ldots, \alpha_k$ is a collection of loops in $A$ that normally generates $L$. Let $X'$ be the complex obtained by attaching a 2-cell along each $\alpha_l$ and let $A'$ be the resulting subcomplex. Then $A' \hookrightarrow X'$ is a homotopy equivalence. (Note: In the above situation, we call $A \hookrightarrow X$ a \textbf{mod L homotopy equivalence}.)
\end{Proposition}

Since $H_*(W''',M; \BZ Q) = 0$, we have that by Proposition \ref{G-T4-lem42}, $W'''$ union the 2-handles $f^2_j$ strong deformation retracts onto $M$ union the 2-handles $f^2_j$. One may now extend via the identity to get a strong deformation retraction $r: W''' \bigcup_{M \bigcup \text{2-handles}} W^{(iv)} \rightarrow W^{(iv)}$. Now $r \circ H$ is the desired strong deformation retraction, of both $W^{(iv)}$ onto its right-hand boundary $N$ and $W^{(iv)}$ read backwards onto its left-hand boundary $N$.

\indexspace

Now, suppose, for the cobordism $(W^{(iv)}, N, M)$, we have $\tau(W^{(iv)},N) = A \ne  0$. As the epimorphism $\eta: G \rightarrow Q$ admist a left inverse $\zeta: Q \rightarrow G$, by the functoriality of Whitehead torsion, we have that $Wh(\eta): Wh(G) \rightarrow Wh(q)$ is onto and admits a left inverse $Wh(\zeta): Wh(q) \rightarrow Wh(G)$. Let $B$ have $A + B = 0$ in $Wh(Q)$ and set $B' = Wh(\zeta)(B)$. By The Realization Theorem from \cite{Rourke-Sanderson}, there is a cobordism $(R, M, N_-)$ with $\tau(R, M) = B'$. If $W = (W^{(iv)} \cup_M R)$, by Theorem 20.2 in \cite{Cohen}, $\tau(W, N) = \tau(W^{(iv)}, N) + \tau(W, W^{(iv)})$. By Theorem 20.3 in \cite{Cohen}, $\tau(W, W^{(iv)}) = Wh(\eta)(\tau(R,M))$. So, $\tau(W^{(iv)}, N) + Wh(\eta)(\tau(R, M) = A +  Wh(\eta)(B') = A + B = 0$, and $(W,N,N_-)$ is a 1-sided s-cobordism.

%

\end{Proof}
        
        \newchapter{USING THE REVERSE PLUS CONSTRUCTION TO BUILD PSEUDO-COLLARS}{USING THE REVERSE PLUS CONSTRUCTION TO BUILD PSEUDO-COLLARS}{USING THE REVERSE PLUS CONSTRUCTION TO BUILD PSEUDO-COLLARS}
    \label{sec:LabelForChapter4}

        
        
    
        \section[Some Preliminaries and the Main Result]{Some Preliminaries and the Main Result}
        \label{sec:LabelForChapter4:Section1}
        Our goal in this section is to display the usefulness of 1-sided s-cobordisms by using them to create large numbers of topologically distinct pseudo-collars (to be defined below), all with similar group-theoretic properties.

\indexspace

We start with some basic definitions and facts concerning pseudo-collars.

\indexspace

\begin{Definition}
Let $W^{n+1}$ be a 1-ended manifold with compact boundary $M^n$. We say $W$ is \textbf{inward tame} if $W$ admits a co-final sequence of ``clean'' neighborhoods of infinity $(N_i)$ such that each $N_i$ is finitely donimated. [A \textbf{neighborhood of infinity} is a subspace the closure of whose complement is compact. A neighborhood of infinity $N$ is \textbf{clean} if (1) $N$ is a closed subset of $W$ (2) $N \cap \partial W = \emptyset$ (3) $N$ is a codimension-0 submanifold with bicollared boundary.]
\end{Definition}

\begin{Definition}
A manifold $N^n$ with compact boundary is a \textbf{homotopy collar} if $\partial N^n \hookrightarrow N^n$ is a homotopy equivalence.
\end{Definition}

\begin{Definition}
A manifold is a \textbf{pseudo-collar} if it is a homotopy collar which contains arbitrarily small homotopy collar neighborhoods of infinity. A manifold is \textbf{pseudo-collarable} if it contains a pseudo-collar neighborhood of infinity.
\end{Definition}

Pseudo-collars naturally break up as 1-sided s-cobordisms. That is, if $N_1 \subseteq N_2$ are homotopy collar neighborhoods of infinity of an end of a pseudo-collarable manifold, the $cl(N_2 \backslash N_1)$ is a cobordism $(W, M, M_-)$, where $M \hookrightarrow W$ is a simple homotopy equivalence. Taking an decreasing chain of homotopy collar neighborhoods of infinity yields a decomposition of a pseudo-collar as a ``stack'' of 1-sided s-cobordisms. 

\indexspace

Conversely, if one starts with a closed manifold $M$ and uses the techniques of chapter 3 to produce a 1-sided s-cobordisms $(W_1, M, M_-)$, then one takes $M_-$ and again uses the techniques of Chapter 3 to produce a 1-sided s-cobordisms $(W_2, M_-, M_{--})$, and so on ad infinitum, and then one glues $W_1 \cup W_2 \cup \ldots$ together to produce an $(n+1)-$dimensional manifold $N^{n+1}$, then $N$ is a pseudo-collar.

\indexspace

So, 1-sided s-cobordisms are the ``correct'' tool to use when constructing pseudo-collars.

\begin{Definition}
The \textbf{fundamental group system at $\infty$}, $\pi_1(\epsilon(X), r)$, of an end $\epsilon(X)$ of a non-compact topological space $X$, is defined by taking a cofinal sequence of neighborhoods of $\infty$ of the end of $X$, $N_1 \supseteq N_2 \supseteq N_3 \supseteq \ldots,$, a proper ray $r: [0, \infty) \rightarrow X$, and looking at its related inverse sequence of fundamental groups $\pi_1(N_1, p_1) \lto \pi_1(N_2, p_2) \lto \pi_1(N_3, p_3) \lto \ldots$ (where the bonding maps are induced by inclusion and the basepoint change isomorphism, induced by the ray $r$).
\end{Definition}

Such a fundamental group system at infinity has a well-defined associated pro-fundamental group system at infinity, given by its equivalence class inside the category of inverse sequences of groups under the below equivalence relation.

\begin{Definition}
Two inverse sequences of groups $(G_i, \alpha_i)$ and $(H_i, \beta_i)$ are said to be \textbf{pro-isomorphic} if there exists subsequences of each, which may be fit into a commuting ladder diagram as follows:
\end{Definition}

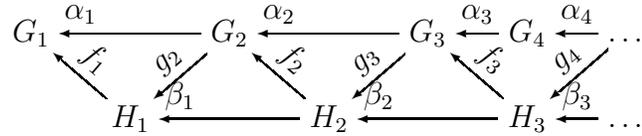
\begin{diagram}[size=14.5pt]
G_1   & \lTo^{\alpha_1} &            &            & G_2        & \lTo^{\alpha_2} &            &            & G_3        & \lTo^{\alpha_3} & G_4        & \lTo^{\alpha_4} & \ldots \\
      & \luTo^{f_1}      &            & \ldTo^{g_2}      &            & \luTo^{f_2}      &            & \ldTo^{g_3}      &            & \luTo^{f_3}      &     & \ldTo^{g_4} & \\
      &            & H_1        & \lTo^{\beta_1} &            &            & H_2        & \lTo^{\beta_2} &            &            & H_3        & \lTo^{\beta_3} & \ldots \\
\end{diagram}

A more detailed introduction to fundamental group systems at infinity can be found in \cite{Geoghegan} or \cite{Guilbault2}.

\begin{Definition}
An inverse sequence of groups is \textbf{stable} if is it pro-isomorphic to a constant sequence $G \lto G \lto G \lto G \ldots$ with the identity for bonding maps.
\end{Definition}

The following is a theorem of Brown from \cite{MBrown}.

\begin{Theorem}
The boundary of a manifold $M$ is collared, i.e., there is a neighborhood $N$ of $\partial M$ in $M$ such that $N \approx \partial M \times \BI$.
\end{Theorem}

The following is from Siebenmann's Thesis, \cite{Siebenmann}.

\begin{Theorem}
An open manifold $W^{n+1}$ ($n \ge 5$) admits a compactification as an $n+1$-dimensional manifold with an $n$-dimensional boundary manifold $M^n$ if \\
(1) $W$ is inward tame \\
(2) $\pi_1(\epsilon(W))$ is stable for each end of $W$, $\epsilon(W)$ \\
(3) $\sigma_{\infty}(\epsilon(W)) \in \widetilde{K}_0[\BZ\pi_1(\epsilon(W))]$ vanishes for each end of $W$, $\epsilon(W)$
\end{Theorem}

\begin{Definition}
An inverse sequence of groups is \textbf{semistable} or \textbf{Mittag-Leffler} if is it pro-isomorphic to a sequence $G_1 \twoheadleftarrow G_2 \twoheadleftarrow G_3 \twoheadleftarrow G_4 \ldots$ with epic bonding maps.
\end{Definition}

\begin{Definition}
An inverse sequence of finitely presented groups is \textbf{perfectly semistable} if and only if is it pro-isomorphic to a sequence $G_1 \twoheadleftarrow G_2 \twoheadleftarrow G_3 \twoheadleftarrow G_4 \ldots$ with epic bonding maps and perfect kernels.
\end{Definition}

The following two lemmas show that optimally chosen perfectly stable inverse sequences behave well under passage to subsequences.

\begin{Lemma}
Let 

\begin{diagram}
1 & \rTo & K & \rTo^{\iota} & G & \rTo^{\sigma} & Q & \rTo &1
\end{diagram}

be a short exact sequence of groups with $K, Q$ perfect. Then $G$ is perfect.
\end{Lemma}

\begin{Proof}
Follows from Lemma 1 in ~\cite{Guilbault1}. Let $g \in G$. Then $\sigma(g) \in Q$, so $\sigma(g) = \Pi_{i=1}^k [x_i,y_i], x_i, y_i \in Q$, as $Q$ is perfect. But, now, $\sigma$ is onto, $\exists u_i \in G$ with $\sigma(u_i) = x_i$ and $v_i \in G$ with $\sigma(v_i) = y_i$. Set $g' = \Pi_{i=1}^k [u_i,v_i]$. Then 

$$\sigma(g \cdot (g')^{-1}) = \sigma(g) \cdot \sigma(g')^{-1} = \Pi_{i=1}^k [x_i,y_i] \cdot (\Pi_{i=1}^k [x_i,y_i])^{-1} = 1 \in Q.$$

 Thus, $g \cdot (g')^{-1} \in \iota(K)$, and $\exists r_j, s_j \in K$ with $g \cdot (g')^{-1} = \iota(\Pi_{j=1}^l [r_j,s_j]$, as $K$ is perfect. But, finally, $g = [g \cdot (g')^{-1}] \cdot g' = \Pi_{j=1}^l [\iota(r_j),\iota(s_j)] \cdot \Pi_{i=1}^k [u_i,v_i]$, which proves $g \in [G,G]$.
\end{Proof}

\begin{Lemma}
If $\alpha: A \rto B$ and $\beta: B \rto C$ are both onto and have perfect kernels, the $(\beta \circ \alpha): A \rto C$ is onto and has perfect kernel.
\end{Lemma}

\begin{Proof}
(Perfect kernel) Set $K = \ker(\alpha), Q = \ker(\beta), G = \ker(\beta \circ \alpha)$

\begin{Claim}
$K = \ker(\alpha|_G): G \rto B$
\end{Claim}

\begin{Proof}
($\subseteq$) Let $g \in G$ have $\alpha(g) = e \in B$ Then $g \in A$ and $\alpha(g) = e \in B$, so $G \in K$

\indexspace

($\supseteq$) Let $k \in K$. Then $\alpha(k) = e \in B$, so $\beta(\alpha(k)) = \beta(e) = e \in Q$. Thus $(\beta \circ \alpha)(k) = e \in C$, and $k \in G$. Since $\alpha(k) = e \in B$, this shows $k \in \ker(\alpha|_G)$.
\end{Proof}

\end{Proof}

The following is a result from \cite{G-T2}.

\begin{Theorem}[Guilbault-Tinsley] \label{G-T2}
A non-compact manifold $W^{n+1}$ with compact (possibly empty) boundary $\partial W = M$ is pseudo-collarable if and only if \\
(1) $W$ is inward tame \\
(2) $\pi_1(\epsilon(W))$ is perfectly semistable for each end of $W$, $\epsilon(W)$ \\
(3) $\sigma_{\infty}(\epsilon(W)) \in \widetilde{K}_0[\BZ\pi_1(\epsilon(W))]$ vanishes for each end of $W$, $\epsilon(W)$
\end{Theorem}

So, the pro-fundamental group system at infinity of a pseudo-collar is perfectly semistable. As is outlined in Chapter 4 of \cite{Guilbault2}, the pro-fundamental group system at infinity is independent of base ray for ends with semistable pro-fundamental group at infinity, and hence for 1-ended pseudo-collars.

\begin{RestateTheorem}{Theorem}{thmpseudo-collars}{Uncountably Many Pseudo-Collars on Closed Manifolds with the Same Boundary and Similar Pro-$\pi_1$}
Let $M^n$ be a closed smooth manifold ($n \ge 6$) with $\pi_1(M) \cong \BZ$ and let $S$ be the fintely presented group $V*V$, which is the free preduct of 2 copies of Thompson's group $V$. Then there exists an uncountable collection of pseudo-collars $\{N^{n+1}_{\omega}\ |\ \omega \in \Omega\}$, no two of which are homeomorphic at infinity, and each of which begins with  $\partial N^{n+1}_{\omega} = M^n$ and is obtained by blowing up countably many times by the same group $S$. In particular, each has fundamental group at infinity that may be represented by an inverse sequence

\begin{diagram}[size=14.5pt]
\BZ & \lOnto^{\alpha_1} & G_{1} & \lOnto^{\alpha_2} & G_{2} & \lOnto^{\alpha_3} & G_{3} & \lOnto^{\alpha_4} & \ldots \\
\end{diagram}

with $ker(\alpha_i) = S$ for all $i$.
\end{RestateTheorem}

We give a brief overview of our strategy. We will start with the manifold $\BS^1 \times \BS^{n-1}$, which has fundamental group $\BZ$. We let $S$ be the free product of 2 copies of Thompson's group $V$, which is a fintely presented, superperfect group for which $Out(S)$ has torsion elements of all orders. Then we will blow $\BZ$ up by $S$ to semi-direct products $G_{p_1}$, $G_{p_2}$, $G_{p_3}$, ..., in infintely many different ways using different outer automorphisms $\phi_{p_i}$ of prime order. We will then use the theorem of last chapter to blow up $\BS^1 \times \BS^{n-1}$ to a manifolds $M_{p_1}$, $M_{p_2}$, $M_{p_3}$, ..., by cobordisms $W_{p_1}$, $W_{p_2}$, $W_{p_3}$, ... . We will then use different automorphisms, each with order a prime number strictly greater than the prime order used in the last step, from the infinite group $Out(S)$ to blow up each of $G_{p_1}$, $G_{p_2}$, $G_{p_3}$, ..., to a different semi-direct products by $S$, and will then use the theorem of last chapter to extend each of $W_{p_1}$, $W_{p_2}$, $W_{p_3}$, ..., in infintely many different ways.
\indexspace

Continuing inductively, we will obtain increasing sequences $\omega$ of prime numbers describing each sequence of 1-sided s-cobordisms. We will then glue together all the semi-s-cobordisms at each stage for each unique increasing sequence of prime numbers $\omega$, creating for each an $(n+1)$-manifold $N^{n+1}_{\omega}$, and show that there are uncountably many such pseudo-collared $(n+1)$-manifolds $N_{\omega}$, one for each increasing sequence of prime numbers $\omega$, all with the same boundary $\BS^1 \times \BS^{n-1}$, and all the result of blowing up $\BZ$ to a semi-direct product by copies of the same superperfect group $S$ at each stage. The fact that no two of these pseudo-collars	are homeomorphic at infinity will follow from the fact that no two of the inverse sequences of groups are pro-isomorphic. Much of the algebra in this chapter is aimed at proving that delicate result.

\begin{Remark}
There is an alternate strategy of blowing up each the fundamental group $G_i$ at each stage by the free product $G_i*S_i$; using a countable collection of freely indecomposible kernel groups $\{S_i\}$ would then allow us to create an uncountable collection of pseudo-collars; an algebraic argument like that found in \cite{Sparks} or \cite {C-K} would complete the proof. However, they would not have the nice kernel properties that our construction has.
\end{Remark}

It seems likely that other groups than Thompson's group $V$ would work for the purpose of creating uncountably many pseudo-collars, all with similar group-theoretic properties, from sequences of 1-sided s-cobordisms. But, for our purposes, $V$ possesses the ideal set of properties.
            
        \section[Some Algebraic Lemmas, Part 1]{Some Algebraic Lemmas, Part 1}
        \label{sec:LabelForChapter4:Section2}
        In this section, we go over the main algebraic lemmas necessary to do our strategy of blowing up the fundamental group at each stage by a semi-direct product with the same superperfect group $S$.

\indexspace

Thompson's group $V$ is finitely presented, superperfect, simple, and contains torsion elements of all orders. Note that simple implies $V$ is centerless, Hopfian, and freely indecomposable.

\indexspace

An introduction to some of the basic properties of Thompson's group $V$ can be found in \cite{CFP}, There, it is shown that $V$ is finitely presented and simple. It is also noted in \cite{CFP} that $V$ contains torsion elements of all orders, as $V$ contains a copy of every symmetric group on $n$ letters, and hence of every finite group. In \cite{Brown2}, it is noted that $V$ is superperfect. We give proofs of some of the simpler properties.

\begin{Lemma}
Every non-Abelian simple group is perfect
\end{Lemma}

\begin{Proof}
Let $G$ be a simple, non-Abelian group, and consider the commutator subgroup $K$ of $G$. This is not the trivial group, as $G$ is non-Abelian, and so by simplicity, must be all of $G$. This shows every element of $G$ can be written as a product of commutator of elements of $G$, and so $G$ is perfect.
\end{Proof}

\begin{Definition}
A group $G$ is \textbf{Hopfian} if every onto map from $G$ to itself is an isomorphism. Equivalently, a group is Hopfian if it is not isomorphic to any of its proper quotients.
\end{Definition}

\begin{Lemma}
Every simple group is Hopfian.
\end{Lemma}

\begin{Proof}
Clearly, the trivial group is Hopfian.

\indexspace

So, let $G$ be a non-trivial simple group. Then the only normal subgroups of $G$ are $G$ itself and $\langle e \rangle$, so the only quotients of $G$ are $\langle e \rangle$ and $G$, respectively. So, the only proper quotient of $G$ is $\langle e \rangle$, which cannot be isomorphic to $G$ as $G$ is nontrivial.
\end{Proof}

Let $S = P_1*P_2$ be the free product of 2 copies of $V$ with itself. This is clearly finitely presented, perfect (by Meyer-Vietoris), and superperfect (again, by Meyer-Vietoris). Note that $S$ is a free product of non-trivial groups, so $S$ is centerless. In \cite{Dey-Neumann}, it is noted that free products of Hofpian, finitely presented, freely indecomposable groups are Hopfian, so $S = V*V$ is Hopfian. $S$ (and not $V$ itself) will be the superperfect group we use in our constructions. 

\indexspace

We need a few lemmas.

\indexspace

\begin{Lemma} \label{lemnontriv-free-prod}
Let $A, B, C, \text{and } D$ be non-trivial groups. Let $\phi: A \times B \rightarrow C*D$ be a surjective homomorphism. Then one of $\phi(A \times \{1\})$ and $\phi(\{1\} \times B)$ is trivial and the other is all of $C*D$
\end{Lemma}

\begin{Proof}

Let $x \in \phi(A \times \{1\}) \cap \phi(\{1\} \times B)$. Then $x \in \phi(A \times \{1\})$, so $x$ commutes with everything in $\phi(\{1\} \times B)$. But $x \in \phi(\{1\} \times B)$, so $x$ commutes with everything in $\phi(A \times \{1\})$. As $\phi$ is onto, this implies $\phi(A \times \{1\}) \cap \phi(\{1\} \times B) \le Z(C*D)$.

\indexspace

But, by a standard normal forms argument, the center of a free product is trivial! So, $\phi(A \times \{1\}) \cap \phi(\{1\} \times B) \le Z(C*D) = 1$. However, this implies that $\phi(A \times \{1\}) \times \phi(\{1\} \times B) = C*D$. By a result in \cite{Baer-Levi}, a non-trivial direct product cannot be a non-trivial free product. (If you'd like to see a proof using the Kurosh Subgroup Theorem, that can be found in many group theory texts, such as Theorem 6.3.10 of \cite{Robinson}.  An alternate, much simpler proof due to P.M. Neumann can be found in \cite{Lyndon-Schupp} in the observation after Lemma IV.1.7). Thus, $\phi(A \times \{1\}) = C*D$ or $\phi(\{1\} \times B) = C*D$ and the other is the trivial group. The result follows.

\end{Proof}

\begin{Corollary}  \label{cornontriv-free-product}
Let $A_1, \ldots, A_n$ be non-trivial groups and let $C*D$ be a free product of non-trivial groups. Let $\phi: A \times \ldots \times A_n \rightarrow C*D$ be a surjective homomorphism.

\indexspace

Then one of the $\phi(\{1\} \times \ldots A_i \times \ldots \times \{1\})$ is all of $C*D$ and the rest are all trivial.

\end{Corollary}

\begin{Proof}

Proof is by induction.

	($n = 2$) This is Lemma \ref{lemnontriv-free-prod}.

(Inductive Step) Suppose the result is true for $n-1$. Set B = $A_1 \times \ldots \times A_{n-1}$. By Lemma \ref{lemnontriv-free-prod}, either $\phi(B \times \{1\})$ is all of $C*D$ and $\phi(\{1\} \times A_n)$ is trivial or $\phi(B \times \{1\})$ is trivial and $\phi(\{1\} \times A_n)$ is all of $C*D$.

\indexspace

If $\phi(B \times \{1\})$ is trivial and $\phi(\{1\} \times A_n)$ is all of $C*D$, we are done.

\indexspace

If  $\phi(B \times \{1\})$ is all of $C*D$ and $\phi(\{1\} \times A_n)$ is trivial, then, by the inductive hypothesis, we are also done.

\end{Proof}

\begin{Corollary} \label{lemstraightening-up-lemma}
Let $S_1, S_2, \ldots, S_n$ all be copies of the same non-trivial free product, and let $\psi: S_1 \times S_2 \times \ldots \times S_n \rightarrow S_1 \times S_2 \times \ldots \times S_n$ be a isomorphism. Then $\psi$ decomposes as a ``matrix of maps'' $\psi_{i,j}$, where each $\psi_{i,j} = \pi_{S_j} \circ \psi|_{S_i}$ (where $\pi_{S_j}$ is projection onto $S_j$), and there is a permutation $\sigma$ on $n$ indices with the property that each $\psi_{\sigma(j), j}: S_{\sigma(j)} \rightarrow S_j$ is an isomorphism, and all other $\psi_{i,j}$'s are the zero map.
\end{Corollary}

\begin{Proof}

By Lemma \ref{cornontriv-free-product} applied to $\pi_{S_j} \circ \psi$, we clearly have a situation where each $\pi_{S_j} \circ \psi|_{S_i}$ is either trivial or onto. If we use a schematic diagram with an arrow from $S_i$ to $S_j$ to indicate non-triviality of a map $\psi_{i,j}$, we obtain a diagram like the following.

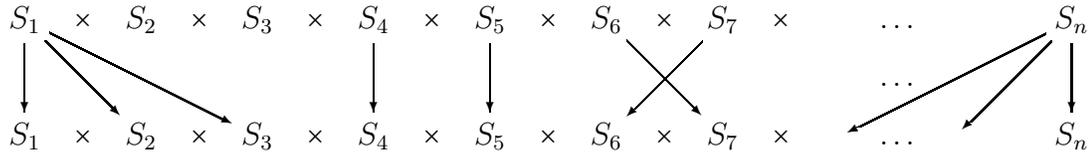
\begin{diagram}[size=22.0pt]
S_1  & \times                & S_2 & \times & S_3  & \times & S_4  & \times & S_5  & \times & S_6 & \times & S_7  & \times & & \ldots &                & & S_n  \\
\dTo & \rdTo(2,2) \rdTo(4,2) &     &        &      &        & \dTo &        & \dTo &        &     & \ldTo(2,2) \rdTo(2,2)    &               &  &            & \ldots & & \ldTo(4,2) \ldTo(2,2) & \dTo \\
S_1  & \times                & S_2 & \times & S_3  & \times & S_4  & \times & S_5  & \times & S_6 & \times & S_7  & \times & & \ldots &                & & S_n  \\
\end{diagram}

where \textit{a priori} some of the $S_i$'s in the domain may map onto multiple $S_j$'s in the target, and there are no arrows eminating from some of the $S_i$'s in the domain.

\indexspace

By the injectivity of $\psi$, there must be at least one arrow eminating from each $S_i$, while by surjectivity of $\psi$, there must be at least one arrow ending at each $S_j$. Corollary \ref{cornontriv-free-product} prevents more than one arrow from ending in a given $S_j$. By the Pidgeonhole Principle, the arrows determine a one-to-one correspondence between the factors in the domain and those in the range. A second application of injectivity now shows each arrow represents an isomorphism.
\end{Proof}

Note that the $\psi_{i,j}$'s form a matrix where each row and each column contain exactly one isomorphism, and the rest of the maps are trivial maps - what would be a permutation matrix (see page 100 in \cite{Robbin}, for instance) if the isomorphisms were replaced by ``1'''s and the trivial maps were replaced by ``0'''s.

\begin{Corollary} \label{corstraightening-up-corollary}
Let $S_1, S_2, \ldots, S_n$ all be copies of the same non-trivial Hopfian free product, and let $\psi: S_1 \times S_2 \times \ldots \times S_n \rightarrow S_1 \times S_2 \times \ldots \times S_m$ be a epimorphism with $m < n$. Then $\psi$ decomposes as a ``matrix of maps'' $\psi_{i,j} = \pi_{S_j} \circ \psi|_{S_i}$, and there is a 1-1 function $\sigma$ from the set $\{1, \ldots, m\}$ to the set $\{1, \ldots, n\}$ with the property that $\psi_{\sigma(j), j}: S_{\sigma(j)} \rightarrow S_{j}$ is an isomorphism, and all other $\psi_{i,j}$'s are the zero map.
\end{Corollary}

\begin{Proof}

Begin with a schematic arrow diagram as we had in the previous lemma. By surjectivity and Lemma \ref{cornontriv-free-product}, each of the $m$ factors in the range is at the end of exactly 1 arrow. From there, we may conclude that each arrow represents an epimorphism, and, hence, by Hopfian, an isomorphism.

\indexspace

To complete the proof, we must argue that at most one arrow can eminate from an $S_i$ factor. Suppose to the contrary, that two arrows emanate from a given $S_i$ factor. Then we have an epimorphism of $S_i$ onto a non-trivial direct product in which each coordinate function is a bijection. This is clearly impossible.
\end{Proof}
        
        \section[Some Algebraic Lemmas, Part 2]{Some Algebraic Lemmas, Part 2}
        \label{sec:LabelForChapter4:Section3}
        Let $\Omega$ be the uncountable set consisting of all increasing sequences of prime numbers $(p_1, p_2, p_3, \ldots)$ with $p_i < p_{i+1}$. For $\omega \in \Omega$ and $n \in \BN$, define $(\omega, n)$ to be the finite sequence consisting of the first $n$ entries of $\omega$.

\indexspace

Let $p_i$ denote the $i^{th}$ prime number, and for the group $S = P_1*P_2$, where each $P_i$ is Thompson's group $V$, choose $u_i \in P_1$ to have $order(u_i) = p_i$.

\indexspace

Recall, if $K$ is a group, $Aut(K)$ is the automorphism group of $K$. Define $\mu: K \rightarrow Aut(K)$ to be $\mu(k)(k') = kk'k^{-1}$. Then the image of $\mu$ in $Aut(K)$ is called \textit{the inner automorphism group of K}, $Inn(K)$. The inner automorphism group of a group $K$ is always normal in $Aut(K)$. The quotient group $Aut(K)/Inn(K)$ is called the \textit{outer automorphism group} $Out(K)$. The kernel of $\mu$ is called the \textit{center of K}, $Z(K)$; it is the set of all $k \in K$ such that for all $k' \in K, kk'k^{-1} = k'$. One has the exact sequence

\indexspace

\begin{diagram}
1 & \rTo & Z(K) & \rTo  & K & \rTo^\mu  & Aut(K) & \rTo^\alpha & Out(K) & \rTo & 1 \\
\end{diagram}

\indexspace

Define a map $\Phi: P_1 \rightarrow Out(P_1*P_2)$ by $\Phi(u) = \phi_u$, where $\phi_u \in Out(P_1*P_2)$ is the outer automorphism defined by the automorphism 

\indexspace

$\phi_{u}(p) = 
\begin{cases}
p & \text{if } p \in P_1 \\
upu^{-1} & \text{if } p \in P_2 \\
\end{cases}$

\indexspace

($\phi_u$ is called a \emph{partial conjugation}.)

\begin{Claim} \label{lemguilbaultembedding-lemma}
$\Phi: P_1 \rightarrow Out(P_1*P_2)$ is an embedding
\end{Claim}

\begin{Proof}
Suppose $\Phi(u)$ is an inner automorphism for some $u$ not $e$ in $P_1$. Since $\Phi(u)$ acts on $P_2$ by conjugation by $u$, to be an inner automorphism, $\Phi(u)$ must also act on $P_1$ by conjugation by $u$. Now, $\Phi(u)$ acts on $P_1$ trivially for all $p \in P_1$, which implies $u$ is in the center of $P_1$. But $P_1$ is centerless! Thus, no $\Phi(u)$ is an inner automophism for any $u \in P_1$.
\end{Proof}

So, for each $u_i$ with prime order the $i^{th}$ prime $p_i$, $\phi_{u_i}$ has prime order $p_i$, as does every conjugate of $\phi_{u_i}$ in $Out(P_1*P_2)$, as $\Phi$ is an embedding.

\begin{Lemma}
For any finite collection of groups $A_1, A_2, \ldots, A_n$, $\Pi_{i=1}^n Out(A_i)$ embeds in $Out(\Pi_{i=1}^n A_1)$.
\end{Lemma}

\begin{Proof}
The natural map from $\Pi_{i=1}^n Aut(A_i)$ to $Aut(\Pi_{i=1}^n A_1)$ which sends a Cartesian product of automorphism individually in each factor to that product considered as an automorphism of the product is clearly an embedding. Now, $Inn(A_1 \times \ldots \times A_n)$ is the image under this natural map of $Inn(A-1) \times \ldots \times Inn(A_n)$, because if $b_i \in A_i$, then $(b_1, \ldots, b_n)^{-1}(a_1, \ldots, a_n)(b_1, \ldots, b_n) = (b_1^{-1}a_1b_1, \ldots, b_n^{-1}a_nb_n)$. So, the induced map on quotient groups, from $\Pi_{i=1}^n Out(A_i)$ to $Out(\Pi_{i=1}^n A_1)$, is also a monomorphism.
\end{Proof}

Now, because the quotient map $\Psi: \Pi_{i=1}^n Out(A_i) \rightarrow Out(\Pi_{i=1}^n A_1)$ is an embedding, $order(\phi_1, \ldots, \phi_n)$ in  $Out(\Pi_{i=1}^n A_1)$ is just $lcm(order(\phi_1), \ldots, order(\phi_n))$, which is just its order in $\Pi_{i=1}^n Out(A_i)$. Moreover each conjugate of $(\phi_1, \ldots, \phi_n)$ in $Out(\Pi_{i=1}^n A_1)$ has the same order $lcm(\phi_1, \ldots, \phi_n)$. Finally, note that if each $\phi_i$ has prime order and each prime occurs only once, then $order(\phi_1, \ldots, \phi_n) = order(\phi_1) \times \ldots \times order(\phi_n)$.

\begin{Lemma} \label{lemguilbaultconjugacy-lemma}
Let $K$ be a group and suppose $\Theta: K \rtimes_{\phi} \BZ  \rightarrow K \rtimes_{\psi} \BZ $ is an isomorphism that restricts to an isomorphism $\overline{\Theta}: K \rightarrow K$. Then $\phi$ and $\psi$ are conjugate as elements of $Out(K)$
\end{Lemma}

\begin{Proof}
We use the presentations $\langle gen(K), a\ |\ rel(K), ak_ia^{-1} = \phi(k_i) \rangle$ and \\ $\langle gen(K), b\ |\ rel(K), bkb^{-1} = \psi(k) \rangle$ of the domain and range respectively, Since $\Theta$ induces an isomorphism on the infinite cyclic quotients by $K$, there exists $c \in K$ with $\Theta(a) = cb^{\pm 1}$. We assume $\Theta(a) = cb$, with the case $\Theta(a) = cb^{-1}$ being similar.

\indexspace

For each $k \in K$, we have

\begin{center}
\begin{tabular}{rcl}
$\Theta(\phi(k))$ & = & $\Theta(aka^{-1})$ \\
                  & = & $\Theta(a)\Theta(k)\Theta(a)^{-1}$ \\
                  & = & $cb\Theta(k)b^{-1}c^{-1}$ \\
                  & = & $c\psi(\Theta(k))c^{-1}$
\end{tabular}
\end{center}

If we let $\iota_c:K \rightarrow K$ denote conjugation by $c$, we have $\overline{\Theta}\phi = \iota_c\psi\overline{\Theta}$ in $Aut(K)$. Quotienting out by $Inn(K)$ and abusing notation slightly, we have $\overline{\Theta}\phi = \psi\overline{\Theta}$ or $\overline{\Theta}\phi\overline{\Theta}^{-1} = \psi$ in $Out(K)$.
\end{Proof}	

\begin{Lemma} \label{lemconderisomorphism-lemma}
For any finite, strictly increasing sequence of primes $(s_1, s_2, \ldots, s_n)$, define $\phi_{(s_1, \ldots, s_n)}: S_1 \times \ldots \times S_n \rightarrow S_1 \times \ldots S_n$ by $\phi_{(s_1, \ldots s_n)}(x_1, \ldots, x_n) = (\phi_{u_1}(x_1), \ldots, \phi_{u_n}(x_n))$, where $\phi_{u_i}$ is the partial conjugation outer automorphism associated above to the element $u_i$ with prime order $s_i$.

\indexspace

Let $(s_1, \ldots, s_n)$ and $(t_i, \ldots, t_n)$ be increasing sequences of prime numbers of length $n$.

\indexspace

Let $G_{(s_1, \ldots, s_n)} = (S_1 \times \ldots \times S_n) \rtimes_{\phi_{(s_1, \ldots, s_n)}} \BZ$ and $G_{(t_i, \ldots, t_n)} = (S_1 \times \ldots \times S_n) \rtimes_{\phi_{(t_1, \ldots, t_n)}} \BZ$ be two semidirect products with such outer actions. Then $G_{(s_1, \ldots, s_n)}$ is isomorphic to $G_{(t_i, \ldots, t_n)}$ if and only if for the underlying sets $\{s_1, \ldots, s_n\} = \{t_1, \ldots, t_n\}$. 
\end{Lemma}

\begin{Proof}
($\Rightarrow$) Let $\theta: G_{(s_1, \ldots, s_n)} \rightarrow G_{(t_i, \ldots, t_n)}$ be an isomorphism. There are $n$ factors of $S$ in the kernel group of each of $G_{(\omega, n)}$ and $G_{(\eta, n)}$. Then $\theta$ must preserve the commutator subgroup, as the commutator subgroup is a characteristic subgroup, and so induces an isomorphism of the perfect kernel group $K = S_1 \times S_2 \times \ldots \times S_n$, say $\overline{\theta}$. By Corollary \ref{lemstraightening-up-lemma}, it must permute the factors of $K$, say via $\sigma$. 

\indexspace

Now, the isomorphism $\theta$ must take the (infinite cyclic) abelianisation \\ $G_{(s_1, \ldots, s_n)}/K_{(s_1, \ldots, s_n)}$ of the one to the (infinite cyclic) abelianisation \\ $G_{(t_i, \ldots, t_n)}/K_{(t_i, \ldots, t_n)}$ of the other, and hence takes a generator of \\ $G_{(s_1, \ldots, s_n)}/K_{(s_1, \ldots, s_n)}$ (say $aK_{(s_1, \ldots, s_n)}$) to a generator of $G_{(t_i, \ldots, t_n)}/K_{(t_i, \ldots, t_n)}$ (say \\ $b^eK_{(t_i, \ldots, t_n)}$, where $bK_{(t_i, \ldots, t_n)}$ is a given generator of $G_{(t_i, \ldots, t_n)}/K_{(t_i, \ldots, t_n)}$ and $e = \pm 1$).  Then since $\theta$ takes $K_{(s_1, \ldots, s_n)} = [G_{(s_1, \ldots, s_n)},G_{(s_1, \ldots, s_n)}]$ to $[G_{(t_i, \ldots, t_n)},G_{(t_i, \ldots, t_n)}] = K_{(t_i, \ldots, t_n)}$, it follows that $\theta$ takes $a$ to a multiple of $b^e$, say $c^{-1}b^e$ where $c$ lies in $K_{(t_i, \ldots, t_n)}$ and $e = \pm 1$. 

\indexspace

Now, by \ref{lemguilbaultconjugacy-lemma}, $\phi_{(s_1, \ldots, s_n)}$ is conjugate in $Out(K)$ to $\phi_{(t_1, \ldots, t_n)}$, $\overline{\theta}(\phi_{(s_1, \ldots, s_n)})\overline{\theta}^{-1} = \phi_{(t_1, \ldots, t_n)}$. But $\Psi$ is an embedding by Lemma \ref{lemguilbaultembedding-lemma}! This shows that $order(\phi_{(s_1, \ldots, s_n)}) = \Pi_{i=1}^n s_i$ and $order(\phi_{(t_1, \ldots, t_n)}) = \Pi_{i=1}^n t_i$ are equal, so, as each $s_i$ and $t_i$ is prime and occurs only once in each increasing sequence, by the Fundamental Theorem of Arithmentic, $\{s_1, \ldots, s_n\} = \{t_1, \ldots, t_n\}$

($\Leftarrow$) Clear.
\end{Proof}

\begin{Lemma} \label{lemconderepimorphism-lemma}
Let $(\omega, n) = (s_1, \ldots, s_n)$ and $(\eta, m) = (t_1, \ldots, t_m)$ be increasing sequences of prime numbers with $n > m$.

\indexspace

Let $G_{(\omega, n)} = (S_1 \times \ldots \times S_n) \rtimes_{\phi_{(\omega, n)}} \BZ$ and $G_{(\eta, m)} = (S_1 \times \ldots \times S_m) \rtimes_{\phi_{(\eta, m)}} \BZ$ be two semidirect products. 	Then there is an epimorphism $g: G_{(\omega, n)} \rightarrow G_{(\eta, m)}$ if and only if $\{t_1, \ldots, t_m\} \subseteq \{s_1, \ldots, s_n\}$. 
\end{Lemma}

\begin{Proof}
The proof in this case is similar to the case $n = m$, except that the epimorphism $g$ must crush out $n-m$ factors of $K_{(\omega,n)} = S_1 \times \ldots \times S_n$ by Corollary \ref{corstraightening-up-corollary} and the Pidgeonhole Principle and then is an isomorphism on the remaining factors.

\indexspace

($\Rightarrow$) Suppose there is an epimorphism $g: G_{(\omega, n)} \rightarrow G_{(\eta, m)}$. Then $g$ must send the commutator subgroup of $G_{(\omega, n)}$ onto the commutator subgroup of $G_{(\eta, m)}$. By Corollary \ref{corstraightening-up-corollary}, $g$ must send $m$ factors of $K_{(\omega,n)} = S_1 \times \ldots \times S_n$ in the domain isomorphically onto the $m$ factors of $K_{(\eta,m)} = S_1 \times \ldots \times S_m$ in the range and sends the remaining $n - m$ factors of $K_{(\omega,n)}$ to the identity. Let $\{i_1, \ldots, i_m\}$ be the indices in $\{1, \ldots, n\}$ of factors in $K_{(\omega,n)}$ which are sent onto a factor in $K_{(\eta,m)} $ and let $\{j_1, \ldots, j_{n-m}\}$ be the indices in $\{1, \ldots, n\}$ of factors in $K_{(\omega,n)}$ which are sent to the identity in $K_{(\eta,m)} $. Then $g$ induces an isomorphism between $S_{i_1} \times \ldots \times S_{i_m}$ and $K_{(\eta,m)} $. Set $L_m = S_{i_1} \times \ldots \times S_{i_m}$

\indexspace

Also, by an argument similar to Lemmas \ref{lemguilbaultconjugacy-lemma} and \ref{lemconderisomorphism-lemma}, $g$ sends sends the infinite cyclic group $G_{(\omega, n)}/K_{(\omega, n)}$ isomorphically onto the infinite cyclic quotient $G_{(\eta, m)}/K_{(\eta, m)}$.

\indexspace

Note that $L_m \rtimes_{\phi_{(s_{i_1}, \ldots, s_{i_m})}} \BZ$ is a quotient group of $G_{(\omega, n)}$ by a quotient map which sends $S_{j_1} \times \ldots \times S_{j_{n-m}}$ to the identity. Consider the induced map $g': L_m \rtimes_{\phi_{(s_{i_1}, \ldots, s_{i_m})}} \BZ \rightarrow G_{(\eta,m)}$. By the facts that $g'$ maps $L_m$ isomorphically onto $K_{(\eta,m)}$ and preserves the infinite cyclic quotients, we have that the kernel of $g$ must equal exactly $S_{j_1} \times \ldots \times S_{j_{n-m}}$; thus, by the First Isomorphism Theorem, we have that $g'$ is an isomorphism.

\indexspace

Finally, $g'$ is an isomorphism of $L_m \rtimes_{\phi_{(s_{i_1}, \ldots, s_{i_m})}} \BZ$ with $G_{(\omega, n)}$ which restricts to an isomorphism of $L_m$ with $S_{t_1} \times \ldots \times S_{t_m}$, so, by Lemma \ref{lemguilbaultconjugacy-lemma}, we have $\phi_{(s_{i_1}, \ldots, s_{i_m})}$ is conjugate to $\phi_{(t_1, \ldots, t_m)}$, so, in $Out(\Pi_{i=1}^n A_1)$, $order(\phi_{(s_{i_1}, \ldots, s_{i_m})}) = order(\phi_{(t_1, \ldots, t_m)})$, and thus, as each $s_i$ and $t_i$ is prime and appears at most once, by an argument similar to Lemma \ref{lemconderisomorphism-lemma} using the Fundamental Theorem of Arithmetic, $\{t_1, \ldots, t_m\} \subseteq \{s_1, \ldots, s_n\}$.

\indexspace

($\Leftarrow$) Suppose $\{t_1, \ldots, t_m\} \subseteq \{s_1, \ldots, s_n\}$. Choose $a \in G_{(\omega, n)}$ with $aK_{(\omega, n)}$ generating the infinite cyclic quotient $G_{(\omega, n)}/K_{(\omega, n)}$ and choose $b \in G_{(\eta, m)}$ with $bK_{(\eta, m)}$ generating the infinite cyclic quotient $G_{(\eta, m)}/K_{(\eta, m)}$. Set $g(a) = b$.

\indexspace

Send each element of $S_i$ (where $S_i$ uses an element of order $t_i$ in its semidirect product definition in the domain) to a corresponding generator of $S_i$ (where $S_i$ uses an element of order $t_i$ in its semidirect product definition in the range) under $g$. Send the elements of all other $S_j$'s to the identity.

\indexspace

Then $g: G_{(\omega, n)} \rightarrow G_{(\eta, m)}$ is an epimorphism. Clearly, $g$ is onto by construction. It remains to show $g$ respects the multiplication in each group.

\indexspace

Clearly, $g$ respects the mutltiplication in each $S_i$ and in $\BZ$

\indexspace

Finally, if $\alpha_i \in S_i$ and $a \in \BZ$,

\begin{tabular}{rcl}
$g(a \alpha_i)$            & = & $g(a)g(\alpha_i)$ \\
$g(\phi_{s_i}(\alpha_i)a)$ & = & $\phi_{t_i}(g(\alpha_i)) g(a)$
\end{tabular}

\indexspace

using the slide relators for each group and the fact that $s_i = t_i$, which implies $\phi_{s_i} = \phi_{t_i}$. So, $g$ respects the multiplication in each group. This completes the proof.
\end{Proof}
        
        \section[Some Algebraic Lemmas, Part 3]{Some Algebraic Lemmas, Part 3}
        \label{sec:LabelForChapter4:Section4}
        Recall $\Omega$ is an uncountable set consisting of increasing sequences of prime numbers $(p_1, p_2, p_3, \ldots)$ with $p_i < p_{i+1}$. For $\omega \in \Omega$ and $n \in \BN$, recall we have defined $(\omega, n)$ to be the finite sequence consisting of the first $n$ entries of $\omega$.

\indexspace

Recall also that $p_i$ denotes the $i^{th}$ prime number, and for the group $S = P_1*P_2$, where each $P_i$ is Thompson's group $V$, we have chosen $u_i \in P_1$ to have $order(u_i) = p_i$.

\indexspace

Recall finally we have define a map $\Phi: P_1 \rightarrow Out(P_1*P_2)$ (where each $P_i$ is a copy of Thompson's group $V$) by $\Phi(u) = \phi_u$, where $\phi_u \in Out(P_1*P_2)$ is the outer automorphism defined by the automorphism 

\indexspace

$\phi_{u}(p) = 
\begin{cases}
p & \text{if } p \in P_1 \\
upu^{-1} & \text{if } p \in P_2 \\
\end{cases}$

\indexspace

(Recall $\phi_u$ is called a \emph{partial conjugation}.)

\indexspace

Set $G_{(\omega, n)} = (S \times S \times  \ldots \times S) \rtimes_{\phi_{(\omega, n)}} \mathbb{Z}$.

\begin{Lemma} \label{lemsessplitting-lemma}
$G_{(\omega, n)} \cong  S \rtimes_{\phi_{w_{s_n}}} G_{(\omega, n-1)}$, where $\phi_{w_{s_n}}$ is partial conjugation by $u_{s_n}$.
\end{Lemma}

\begin{Proof}
First, note that there is a short exact sequence 

\begin{diagram}
1 & \rTo & S & \rTo^{\iota} & G_{(\omega, n)} & \rTo^{\alpha_{n}} & G_{(\omega, n-1)} & \rTo & 1
\end{diagram}

where $\iota$ takes $S$ identically onto the $n^{\text{th}}$ factor, and $\alpha$ crushes out factor, as described in Lemma \ref{lemconderepimorphism-lemma}.

\indexspace

Next, note that there is a left inverse  $j: G_{(\omega, n-1)} \rightarrow G_{(\omega, n)}$ to $\alpha$ given by (1) sending the generator $a$ of the $\BZ$ from its image $\gamma_{n-1}(a)$ in the semi-direct product 

\begin{diagram}
1 & \rTo & (S \times \ldots \times S) & \rTo^{\iota_{n-1}} & G_{(\omega, n-1)} & \rTo^{\beta_{n-1}} & \BZ & \rTo & 1
\end{diagram}

where $\gamma_{n-1}$ is a left inverse to $\beta_{n-1}$, to its image $\gamma_n(a)$ in 

\begin{diagram}
1 & \rTo & (S \times \ldots \times S) & \rTo^{\iota_{n}} & G_{(\omega, n)} & \rTo^{\beta_{n}} & \BZ & \rTo & 1
\end{diagram}

where $\gamma_{n}$ is a left inverse to $\beta_{n}$

and (2) sending each of the images $\iota_{n-1}(t_i)$ of the elements $t_i$ of the $S_i$ associated with $\phi_{w_{s_i}}$ in $G_{(\omega, n-1)}$ in 

\begin{diagram}
1 & \rTo & (S \times \ldots \times S) & \rTo^{\iota_{n-1}} & G_{(\omega, n-1)} & \rTo^{\beta_{n-1}} & \BZ & \rTo & 1
\end{diagram}

to the images $\iota_n(t_i)$ of the elements	 $t_i$ of the $S_i$ associated with $\phi_{s_i}$ in $G_{(\omega, n)}$ in

\begin{diagram}
1 & \rTo & (S \times \ldots \times S) & \rTo^{\iota_{n}} & G_{(\omega, n)} & \rTo^{\beta_{n}} & \BZ & \rTo & 1
\end{diagram}

for $i \in \{1, \ldots, n-1\}$

\indexspace

The existence of a left inverse proves the group extension is a semi-direct product. The needed outer action for the final copy of $S$ in $G_{(\omega, n)}$ may now be read off the defining data for $G_{(\omega, n)}$ in the definition $G_{(\omega, n)} = (S \times S \times  \ldots \times S) \rtimes_{\phi_{(\omega, n)}} \mathbb{Z}$, showing that it is indeed partial conjugation by $u_{s_n}$.

\indexspace

(Alternately, one may note there is a presentation for $(S \times S \times  \ldots \times S) \rtimes_{\phi_{(\omega, n)}} \mathbb{Z}$ that contains a presentation for $S \rtimes_{\phi_{w_{s_n}}} G_{(\omega, n-1)}$

\indexspace

Generators: $z$, the generator of $\BZ$, together with the generators of the first copy of $S$, the generators of the second copy of $S$, ..., and the generators of the $n^{\text{th}}$ copy of $S$.

\indexspace

Relators defining $P_i$'s: the relators for the copy of $P_1$ in the first copy of $S$, the relators for the copy of $P_2$ in the first copy of $S$, the relators for the copy of $P_1$ in the second copy of S, the relators for the copy of $P_2$ in the second copy of $S$, , ..., and the relators for the copy of $P_1$ in the $n^{\text{th}}$ copy of $S$, the relators for the copy of $P_2$ in the $n^{\text{th}}$ copy of $S$.

\indexspace 

Slide Relators: The slide relators between $z$ and the generators of $P_2$ in the first copy of $S$ due to the semi-direct product, the slide relators between $z$ and the generators of $P_2$ in the second copy of $S$ due to the semi-direct product, ..., the slide relators between $z$ and the generators of $P_1$ in the $n^{\text{th}}$ copy of $S$ due to the semi-direct product, and the slide relators between $z$ and the generators of $P_2$ in the $n^{\text{th}}$ copy of $S$ due to the semi-direct product.)
\end{Proof}

Now, this way of looking at $G_{(\omega, n)}$ as a semi-direct product of $S$ with $G_{(\omega, n-1)}$ yields an inverse sequence $(G_{(\omega, n)}, \alpha_n)$, which looks like

\begin{diagram}
G_{(\omega, 0)} & \lTo^{\alpha_{0}} & G_{(\omega, 1)} & \lTo^{\alpha_{1}} & G_{(\omega, 2)} & \lTo^{\alpha_{2}} & \ldots 
\end{diagram}

with bonding maps $\alpha_i: G_{(\omega, i+1)} \rightarrow G_{(\omega, i)}$ that each crush out the most recently added copy of $S$.

\indexspace

A subsequence will look like 

\begin{diagram}
G_{(\omega, n_0)} & \lTo^{\alpha_{n_0}} & G_{(\omega, n_1)} & \lTo^{\alpha_{n_1}} & G_{(\omega, n_2)} & \lTo^{\alpha_{n_2}} & \ldots 
\end{diagram}

with bonding maps $\alpha_{n_i}: G_{\omega, n_j)} \rightarrow G_{\omega, n_i)}$ that each crush out the most recently added $n_j - n_i$ copies of $S$.

\indexspace

\begin{Lemma} \label{lemladderdiagramsnonequiv-lemma}
If, for inverse sequences $(G_{(\omega, n)}, \alpha_{n})$, where $\alpha_{n}: G_{(\omega, n)} \rightarrow G_{(\omega, n-1)}$ is the bonding map crushing out the most recently-added copy of $S$, $\omega$ does not equal $\eta$, then the two inverse sequences are not pro-isomorphic.
\end{Lemma}

\indexspace

\begin{Proof}
Let $(G_{(\omega, n)}, \alpha_{n})$ and $(G_{(\eta, m)}, \beta_{m})$ be two inverse sequences of group extensions, assume there exists a commuting ladder diagram between subsequences of the two, as shown below. By discarding some terms if necessary, arrage that $\omega$ and $\eta$ do not agree beyond the term $n_0$.

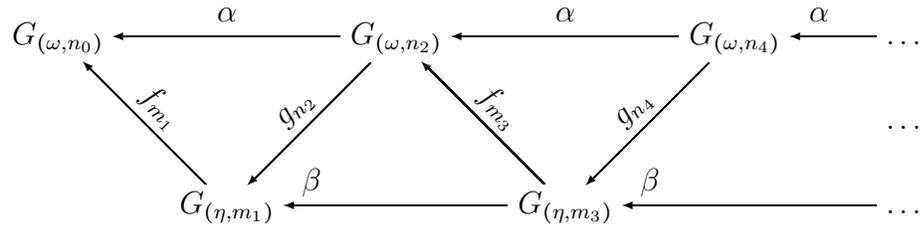
\begin{diagram}[size=32pt]
G_{(\omega, n_0)} &         & \lTo^{\alpha}    &              & G_{(\omega, n_2)} &         & \lTo^{\alpha}    &              & G_{(\omega, n_4)} & \lTo^{\alpha} & \ldots \\
                  & \luTo^{f_{m_1}} &                  & \ldTo^{g_{n_2}}      &                   & \luTo^{f_{m_3}} &                  & \ldTo^{g_{n_4}}      &                   &               & \ldots \\
                  &         &  G_{(\eta, m_1)} & \lTo^{\beta} &                   &         &  G_{(\eta, m_3)} & \lTo^{\beta} &                   &               & \ldots \\
\end{diagram}

By the commutativity of the diagram, all $f$'s and $g$'s must be epimorphisms, as all the $\alpha$'s and $\beta$'s are.

\indexspace

Now, it is possible that $g_{(\omega, n_2)}$ is an epimorphism; by Lemma \ref{lemconderepimorphism-lemma}, $(\eta, m_1)$ might be a subset of $(\omega, n_2)$ when considered as sets. But, $f_{(\eta, m_3)}$ cannot also be an epimorphism, since $(\omega, n_2)$ cannot be a subset of $(\eta, m_3)$ when considered as sets. Since the two sequences can only agree up to $n_0$, if $(\eta, m_1)$ is a subset of $(\omega, n_2)$ when considered as sets, then there must be an prime $p_i$ in $(\omega, n_2)$ in between some of the primes of $(\eta, m_1)$. This prime $p_i$ now cannot be in $(\eta, m_3)$ and is in $(\omega, n_2)$, so we cannot have $(\omega, n_2)$ a subset of $(\eta, m_3)$ when considered as sets, so $f_{(\eta, m_3)}$ cannot be an epimorphism. 

%
%
%
%
%
\end{Proof}
        
        \section[Manifold Topology]{Manifold Topology}
        \label{sec:LabelForChapter4:Section5}
        We now begin an exposition of our example.

\indexspace

\begin{RestateTheorem}{Theorem}{thmpseudo-collars}{Uncountably Many Pseudo-Collars on Closed Manifolds with the Same Boundary and Similar Pro-$\pi_1$}
Let $M^n$ be a closed smooth manifold ($n \ge 6$) with $\pi_1(M) \cong \BZ$ and let $S$ be the fintely presented group $V*V$, which is the free preduct of 2 copies of Thompson's group $V$. Then there exists an uncountable collection of pseudo-collars $\{N^{n+1}_{\omega}\ |\ \omega \in \Omega\}$, no two of which are homeomorphic at infinity, and each of which begins with  $\partial N^{n+1}_{\omega} = M^n$ and is obtained by blowing up countably many times by the same group $S$. In particular, each has fundamental group at infinity that may be represented by an inverse sequence

\begin{diagram}[size=14.5pt]
\BZ & \lOnto^{\alpha_1} & G_{1} & \lOnto^{\alpha_2} & G_{2} & \lOnto^{\alpha_3} & G_{3} & \lOnto^{\alpha_4} & \ldots \\
\end{diagram}

with $ker(\alpha_i) = S$ for all $i$.
\end{RestateTheorem}

\begin{Proof}
For each element $\omega \in \Omega$, the set of all increasing sequences of prime numbers, we will construct a pseudo-collar $N_{\omega}^{n+1}$ whose fundamental group at infinity is represented by the inverse sequence $(G_{(\omega, n)}, \alpha_{(\omega, n)})$. By Lemma \ref{lemladderdiagramsnonequiv-lemma}, no two of these pseudo-collars can be homeomorphic at infinity, and the Theorem will follow.

\indexspace

To form one of the pseudo-collars, start with $M = \BS^1 \times \BS^{n-1}$ with fundamental group $\BZ$ and then blow it up, using Theorem \ref{thmsemi-s-cob}, to a cobordism $(W_{(s_1)}, M, M_{(s_1)})$ corresponding to the group $G_{(s_1)}$ ($s_1$ a prime)..

\indexspace

We then blow this right-hand boundaries up, again using Theorem \ref{thmsemi-s-cob} and Lemma \ref{lemsessplitting-lemma}, to cobordisms $(W_{(s_1, s_2)}, M_{(s_1)}, M_{(s_1, s_2)})$ corresponding to the group $G_{(s_1, s_2)}$ above.

\indexspace

We continue in the fashion \textit{ad infinitum}.

\indexspace

The structure of the collection of all pseudo-collars will be the set $\Omega$ described above.

\indexspace

We have shown that the pro-fundamental group systems at infinity of each pseudo-collar are non-pro-isomorphic in Lemma \ref{lemladderdiagramsnonequiv-lemma}, so that all the ends are non-diffeomorphic (indeed, non-homeomorphic).

\indexspace

This proves we have uncountably many pseudo-collars, each with boundary $M$, which have distinct ends.
\end{Proof}

\begin{Remark}
The above argument should generalize to any manifold $M^n$ with $n \ge 6$ where $\pi_1(M)$ is a finitely generated Abelian group of rank at least 1 and any finitely presented, superperfect, centerless, freely indecomposable, Hopfian group $P$ with an infinite list of elements of different orders (the orders all being prime numbers was a convenient but inessential hypothesis). The quotient needs to be Abelian so that the commutator subgroup will be the kernel group, which is necessarily superperfect; the quotient group must have rank at least 1 so that there is an element to send into the kernel group to act via the partial conjugation. The rest of the conditions should be self-explanatory.
\end{Remark}

    %
    %

    \newchapter{EXTENSIONS}{SOME ONE-ENDED MANIFOLDS WHICH ARE NOT PSEUDO-COLLARABLE}{SOME ONE-ENDED MANIFOLDS WHICH ARE NOT PSEUDO-COLLARABLE}
    \label{sec:LabelForChapter5}

        
        
    
        \section[Hypo-Abelian Groups]{Hypo-Abelian Groups}
        \label{sec:LabelForChapter5:Section1}
        In \cite{G-T2}, Guilbault and Tinsley construct the first known example of an inward tame but non-pseudo-collarable 1-ended manifolds with compact boundary. In fact, their example satifies conditions (1) and (3) of Theorem \ref{G-T2} (a condition hereinafter referred to as \emph{absolutely inward tame}), but fails to be pseudocollarable because it does not satisfy condition (2) of Theorem \ref{G-T2}, that is, it does not have perfectly semistable pro-fundamental group at infinity. Their example is based on a single inverse sequence of groups, created specifically for the purpose of their example.

\indexspace

In this chapter, we present a more general strategy for creating absolutely inward tame manifolds that are not pseudocollarable. Our construction begins with any manifold $M^n$, $n \ge 5$, whose fundamental group is hypo-Abelian (to be defined below) and which contains an element of infinite order. From this, we create a homotopy collar $W^{n+1}$ whose boundary $\partial W = M$ and is absolutely inward tame but is not pseudo-collarable.

\begin{Definition}
A group G is said to \textbf{hypo-Abelian} if it contains no non-trivial perfect subgroup.
\end{Definition}

\indexspace

Examples:

(1) If $G$ is Abelian, then $G$ is hypo-Abelian.

\indexspace

(2) Recall that a group is \textit{solvable} if and only if its derived series

$$G = G_0 \trianglerighteq G^{(1)} \trianglerighteq G^{(2)} \trianglerighteq \ldots \trianglerighteq G^{(n)} = \langle e \rangle$$

terminates at a finite length in the trivial group, where each $G^{(i)} = [G^{(i-1)}, G^{(i-1)}]$.

\indexspace

Since the (possibly transfinite) derived series of a group always terminates in the perfect core of the group (the largest perfect subgroup of the group), and solvable groups have their derived series terminate in the trivial group, every solvable group is hypo-Abelian.

\indexspace

(3) Free groups are hypo-Abelian

\indexspace

Since by the Nielsen-Schreier Theorem (see \cite{Neilsen} and \cite{Schreier}), subgroups of free groups are free.

\indexspace

(4) Free products of hypo-Abelian groups are hypo-Abelian.

\indexspace

This is because of the Kurosh Subgroup Theorem, which states that a subgroup of a free product is $F *_{\lambda \in \Lambda} \alpha_{\lambda}^{-1}P_{\lambda}\alpha_{\lambda}$, where $F$ is a free group and $\alpha_{\lambda}^{-1}P_{\lambda}\alpha_{\lambda}$ is a conjugate of a subgroup of one of the groups in the free product. Since free groups admit no non-trivial perfect subgroups and each factor admits no non-trivial perfect subgroup, there are no non-trivial perfect subgroups, and free products of hypo-Abelian groups are hypo-Abelian.

\indexspace

(5) Every extension of a hypo-Abelian group by a hypo-Abelian group is hypo-Abelian.

\begin{Lemma} \label{hypo-ext}
Let 

\begin{diagram}
1 & \rTo & K & \rTo^{\iota} & G & \rTo^{\sigma} & Q & \rTo &1
\end{diagram}

be a short exact sequence of groups with $K, Q$ hypo-Abelian. Then $G$ is hypo-Abelian.
\end{Lemma}

\begin{Proof}

Let $P \le G$ be perfect. Then $\sigma(P)$ is perfect, for if $x = [y,z]$ in $G$, then $\sigma(x) = [\sigma(y),\sigma(z)]$ in $Q$. But $Q$ is hypo-Abelian, so $\sigma(P) =\ \langle e \rangle\ \le\ Q$. This means $P \le \iota(K)$. But $K$ is hypo-Abelian, so $\iota(K)$ is hypo-Abelian, so $P =\ \langle e \rangle$.

\end{Proof}

\indexspace

(6) In \cite{Howie2}, Howie shows that every right-angled Artin group (RAAG) is hypo-Abelian.

\indexspace

(7) Split amalgamated free products are hypo-Abelian

\begin{Definition}
A monomorphism $\alpha: A \hookrightarrow B$ is said to be \textbf{split} if there is a homomorphism $\beta: B \twoheadrightarrow A$ such that $\beta \circ \alpha = id_A$.
\end{Definition}

\begin{Definition}
A group $G$ is a \textbf{split amalgamated free product} if and only if $G$ may be expressed as an amalgamated free product $B *_{A} C$, where one of the injections $A \hookrightarrow B$ or $A \hookrightarrow C$ splits.
\end{Definition}

The following is a theorem from \cite{Howie1}.

\begin{Theorem}[Theorem E] \label{thmhowie-thm-e}
Let $T$ denote the class of hypo-Abelian groups. Then $T$ is closed under the operation of split amalgamated free products.
\end{Theorem}

(The above result will be important in proving the main theorem of this chapter.)

(8) Every residually solvable group is hypo-Abelian.

\indexspace

(9) Every group $G$ has a ``hypo-Abelianization'' obtained by quotienting out its perfect core, the largest perfect subgroup.

\indexspace

(10) The Baumslag-Solitar groups BS(1,n) given by $\langle t, x\ |\ x = tx^nt^{-1} \rangle$.

\indexspace

The Baumslag-Solitar group BS(1,n), $I_n$, fits into a short exact sequence

\begin{diagram}
1 & \rTo & \BZ[\frac{1}{n}] & \rTo & I_n & \rTo^{r} & \BZ & \rTo & 1
\end{diagram}
 
To see this, note that a K(G,1) for $I_n$ is given by taking a cylinder $S^1 \times [0,1]$ and glueing one end to the other by a degree $n$ map; the covering space corresponding to unraveling the generator going across the cylinder is a bi-infinite mapping telescope of degree $n$ maps of the circle to itself, which is easily seen to have the fundamental group isomomorphic to $\BZ[\frac{1}{n}]$. The short exact sequence follows. By Lemma \ref{hypo-ext}, $I_n$ is now hypo-Abelian.
        
        
        \section[Algebriac Lemmas]{Algebraic Lemmas}
        \label{sec:LabelForChapter5:Section2}
        The following establishes the group theory basis for our geometric example of a 1-ended manifold $W$ ``built'' from a closed manifold $M$ with hypo-Abelian fundamental group $G_0$ and an element of infinite order $t_0$. We construct inductively an inverse sequence of amalgamated free products, first of $G_0$ with the Baumslag-Solitar group BS(1,2), $I$, and then of the newly constructed group $G_j$ again with the Baumslag-Solitar group BS(1,2), $I_2$. This inverse sequence of groups will be the fundamental group at infinity of our manifold $W$.

\begin{Theorem} \label{group-theory}
Let $G_0 = \langle A_0\ |\ R_0 \rangle$ be a finitely presented hypo-Abelian group that contains an element $t_0$ of infinite order. Then there is an inverse sequence of hypo-Abelian groups $G_0 \twoheadleftarrow G_1 \twoheadleftarrow G_2 \twoheadleftarrow G_2 \twoheadleftarrow \ldots$ such that for each $i > 0$, $G_j = \langle A_{j-1}, t_j\ |\ R_{j-1}, t_j = [t_j, t_{j-1}] \rangle$. 
\end{Theorem}

\begin{Proof} First note that $G_1 =\ \langle A_{0}, t_1\ |\ R_{0}, t_1 = t_1^{-1}t_{0}^{-1}t_1t_{0} \rangle\ =\ \langle A_{0}, t_1\ |\ R_{0}, t_1^2 = t_{0}^{-1}t_1t_{0} \rangle\ =\ \langle A_{0}, t_1\ |\ R_{0}, t_1 = t_{0}t_1^2t_{0} \rangle \  \cong G_{0} *_{\langle t_{0}\rangle} I_2$, where $I_2$ is the Baumslag-Solitar group BS(1,2) given by $\langle t, x\ |\ x = tx^2t^{-1} \rangle$, with the generator $t \in I_2$ identified with $t_0 \in G_0$, $x \in I_2$ identified with $t_{1} \in G_{1}$. We must show $t_0$ has infinite order in $G_1 \cong G_{0} *_{\langle t_{0}\rangle} I_2$ and $t_1$ has infinite order in $G_1$.

\indexspace

The fact that $t_0$ has infinite order in $G_1$ follows from the facts that $t_0$ has infinite order in $G_0$, $t$ has infinite order in $I_2$, and Britton's Lemma, which states that each factor in a free product with amalgamation embeds in the free product with amalgamation.

\indexspace

Similarly, the fact that $t_1$ has infinite order in $G_1$ follows from the fact that $x$ has infinite order in $I_2$ and Britton's Lemma.

\indexspace

Note that $G_j =\ \langle A_{j-1}, t_j\ |\ R_{j-1}, t_j = t_j^{-1}t_{j-1}^{-1}t_jt_{j-1} \rangle\ =\ \langle A_{j-1}, t_j\ |\ R_{j-1}, t_j^2 = t_{j-1}^{-1}t_jt_{j-1} \rangle\ =\ \langle A_{j-1}, t_j\ |\ R_{j-1}, t_j = t_{j-1}t_j^2t_{j-1} \rangle \  \cong G_{j-1} *_{\langle t_{j-1}\rangle} I$, where $I$ is the Baumslag-Solitar group BS(1,2) given by $\langle t, x\ |\ x = tx^2t^{-1} \rangle$, with the generator $t \in I$ identified with $t_0 \in G_0$, $x \in I$ identified with $t_{1} \in G_{1}$, and where $t_{j}$ has infinite order in $G_{j}$. Note that $t_{j-1}$ has infinite order in $G_j$ as $t_{j-1}$ has infinte order in $G_{j-1}$ by the inductive hypothesis, $t$ has infinite order in $I_2$, and Britton's Lemma again.

\indexspace

The onto map $r: I \twoheadrightarrow \BZ$ above induces an onto map $r_j: G_j \twoheadrightarrow G_{j-1}$ of the free product with amalgamation $G_j \cong G_{j-1} *_{\langle t_{j-1}\rangle} I$. 

\indexspace

By induction, Theorem E from \cite{Howie1} tells us that each $G_j$ is hypo-Abelian.

\indexspace

It remains to show that $t_j$ has infinite order in $G_j$. This follows from the fact that $x$ has infinite order in $I$ (again, by considering the earlier short exact sequence for $I$) and Britton's Lemma.

\end{Proof}

\begin{Remark}
The inverse sequence used by Guilbault-Tinsley in \cite{G-T2} may now be viewed as a special case of Theorem \ref{group-theory} in which $G_0$ is infinite cyclic.
\end{Remark}

We have a very geometric proof that all the $G_j$'s are hypo-Abelian, derived before we knew of Howie's theorem. We give a brief sketch of the proof here.

\begin{Fact} \label{fact-munkres}
Let $p: \hat{X}(H) \rightarrow X$ be a regular covering projection of CW complexes, $Z \subseteq X$ be a connected subcomplex. Then $p^{-1}(Z) \rightarrow Z$ is a covering map 
\end{Fact}

\begin{Proof}

This is Theorem 53.2 in \cite{Munkres}.

\end{Proof}

Let $X_{j-1}$ be a $K(G_{j-1},1)$. Let $p: \widetilde{X_{j-1}} \rightarrow X_{j-1}$ be the universal cover of $X_{j-1}$. Let $Z$ be a simple closed curve in $X_{j-1}$ representing $t_{j-1}$. Let $A_{j-1} = \{t_{j-1}; G_{j-1}\} \cong \mathbb{Z}$. Let $H_{j-1} = \{gA_{j-1}\ |\ g \in G_{j-1}/\BZ\}$ be an indexing set for the collection of path components of $p^{-1}(Z)$. Then $p^{-1}(Z) \equiv \{h\mathbb{R}\ |\ h \in H_{j-1}\}$ as $t_{j-1}$ has infinite order in $G_{j-1}$ and by Fact \ref{fact-munkres} is a covering map.

\indexspace

Let $Y$ be the $K(I,1)$ described earlier. Let $q: \hat{Y}(\BZ[\frac{1}{2}]) \rightarrow Y$ be the intermediate cover of $Y$ with corresponding to the normal subgroup $\BZ[\frac{1}{2}]$ of $I$. Let $Z'$ be a loop in $Y$ representing $t_{j}$, a generator going across the cylinder. Then $q^{-1}(Z')$ is a single line $\mathbb{R}$ in $\hat{Y}(\BZ[\frac{1}{2}])$ as $t_{j}$ has infinite order in $I$ and by Fact \ref{fact-munkres} is a covering map.

\indexspace

Attach a copy of $\hat{Y}(\BZ[\frac{1}{2}])$ (briefly, $\hat{Y}$) along its copy of $\mathbb{R}$ to each component of $p^{-1}(Z)$ in $\widetilde{X_{j-1}}$, taking care to match up basepoint with copy of basepoint and so that the image of 0 under the deck transformation taking $0 \in \mathbb{R}$ to $1 \in \mathbb{R}$ in the copy of $\mathbb{R}$ in $p^{-1}(Z)$ in $\widetilde{X_{j-1}}$ matches up with the image of 0 under the deck transformation taking $0 \in \mathbb{R}$ to $1 \in \mathbb{R}$ in the copy of $\mathbb{R}$ in the copy of $\hat{Y}(\BZ[\frac{1}{2}])$. Call the resulting space $Q$. Let $X_j$ be a the adjunction space formed by gluing $X_{j-1}$ to $Y$ along $Z^{-2}Z'^{-1}ZZ'$, where $Z'$ is the loop in $Y$ reprsesenting $t_j$ and $Z$ is the loop in $X_{j-1}$ representing $t_{j-1}$.

\indexspace

An elementary but tedious argument gives the following.

\begin{Claim}
(1) $q: Q \rightarrow X_j$ evenly covers $X_j$ (with an appropriate adjunction map $q$ as covering map) \\
(2) $q: Q \rightarrow X_j$ is a regular cover \\
(3) There is a homomorphism $\Psi: \pi_1(Q,\widetilde{*}) \cong \ker(r_j)$ and an isomorphism $\Phi:$ Deck($r$) $\rightarrow G_{j-1}$, $G_{j-1} = \pi_1(X_{j-1})$, which makes the following diagram commute:

\begin{diagram}[size=14.5pt]
1 & \rTo & \pi_1(Q, \widetilde{*}) & \rTo^{\beta} & \pi_1(X_j, *) & \rTo^{\alpha} & \text{Deck(}r\text{)} & \rTo & 1 \\
  &      & \dTo^{\Psi}             &              & \dTo^{\Lambda} &               & \dTo^{\Phi}                    &      &   \\
1 & \rTo & \ker(r_j)            & \rTo^{\iota} & G_j           & \rTo^{r_j} & G_{j-1}                  & \rTo & 1 \\
\end{diagram}

where $\alpha$, $\beta$, and $\Lambda$ are the canonical maps from covering space theory. 
\end{Claim}

Since $\widetilde{X}_{j-1}$ is contractible, it is easy to check that $\pi_1(Q)$ is a possibly infinite free product of copies of the Abelian group $\pi_1(\hat{Y}(K)) \cong \BZ[\frac{1}{2}]$, which, by parts (1) and (4) of the examples at the beginning of the section, is hypo-Abelian. We may now apply Lemma \ref{hypo-ext} to conclude $G_j$ is hypo-Abelian.

\indexspace

Recall Lemma 4.1 from \cite{G-T2}.

\begin{Lemma} \label{lemhypo-not-per-semi}
Let $G_0 \twoheadleftarrow G_1 \twoheadleftarrow G_2 \twoheadleftarrow G_3 \twoheadleftarrow G_4 \ldots$ be an inverse sequence of hypo-Abelian groups and non-injective epimorphisms. Then this inverse sequence is not perfectly semistable. 
\end{Lemma}
        
        \section[Manifold Topology]{Manifold Topology}
        \label{sec:LabelForChapter5:Section3}
        Our primary contribution to the following theorem is contained in the algebra presented above. Since the handle-theoretic construction is nearly identical to that provided in Theorem 4.4 in \cite{G-T2}, we provide only an outline. The reader is referred to \cite{G-T2} for details.

\begin{RestateTheorem}{Theorem}{thm-non-pcm}{Existence of Non-Pseudo-Collarable ``Nice'' Manifolds}
Let $M^n$ be an orientable, closed manifold ($n \ge 6$) such that $\pi_1(M)$ contains an element $t_0$ of infinite order and $\pi_1(M)$ is hypo-Abelian. Then there exists a 1-ended, orientable manifold $W^{n+1}$ with $\partial W = M$ in which all clean neighborhoods of infinity have finite homotopy type, but which does not have perfectly semistable fundamental group at infinity. Thus, $W^{n+1}$ is absolutely inward tame but not pseudocollable.
\end{RestateTheorem}

\begin{Sketch}

We need to construct $W^{n+1}$. We will construct $W^{n+1}$ as promised so that a representative  of pro-$\pi_1$ is the inverse sequence  provided by Theorem \ref{group-theory} with $\pi_1(M)$ playing the role of $G_0$. Then Lemma \ref{lemhypo-not-per-semi} will guarantee that $W$ does not have perfectly semi-stable fundamental group at infinity.

\indexspace

We will construct $W^{n+1}$ by creating a sequence of compact cobordisms $(W_i, M_{i}, M_{i+1})$ such that

\begin{enumerate}[a)]
	\item The left-hand boundary of $W_0$ is $M$ with $\pi_1(M)$ hypo-Abelian and an element $t_0$ of infinite order, and, for all $i \ge 1$, the left-hand boundary of $W_i$, $M_{i}$, is the same as the right-hand boundary of $W_{i-1}$.
	\item For all $i \ge 0$, $\pi_1(M_i) \cong G_i$, and $M_i \hookrightarrow W_i$ induces a $\pi_1$ isomorphism.
	\item The isomorphism between $\pi_1(M_i)$ abd $G_i$ can be chosen so that we have a commutative diagram
\begin{diagram}[size=14.5pt]
G_{i-1}        & \lTo & G_i           \\
\dTo^{\cong}    &      & \dTo^{\cong} \\
\pi_1(M_{i-1}) & \lTo & \pi_1(M_i)    \\
\end{diagram}
\end{enumerate}

We will let 

$$W^{n+1} = W_0 \cup W_1 \cup W_2 \cup W_3 \cup \ldots$$

Then for each $i \ge 1$, 

$$N_i = W_i \cup W_{i+1} \cup W_{i+2} \cup \ldots$$

is a clean, connected neighborhood of infinity. Moreover, by properties b) and c) and repeated applications of Seifert-Van Kampen, the inverse system 

$$\pi_1(N_1, p_1) \leftarrow \pi_1(N_2, p_2) \leftarrow \pi_1(N_3, p_3) \leftarrow \ldots$$

is isomorphic to the inverse sequence from Theorem \ref{group-theory}.

\indexspace

(i = 1) Start with $M^n$ and cross it with $\mathbb{I}$. Attach a trivial 1-handle $\alpha_1^1$ corresponding to $t_1$ to the right-hand boundary (abbreviated RHB) of $M \times \mathbb{I}$. Let $t_0$ be the element of infinite order in $\pi_1(M)$. Attach a 2-handle $\alpha_1^2$ for the relator $t_1 = [t_1, t_0]$ to the RHB of $M \times \mathbb{I}$. Set $B_1^{n+1}$ to be $(M \times \mathbb{I}) \cup \alpha_1^1 \cup \alpha_1^2$. 

\begin{Claim}
Then $M_1 \hookrightarrow B_1$ induces a $\pi_1$ isomorphism.
\end{Claim}

\begin{Proof}

By inverting the handlebody decomposition, we may view $B_1$ as the result of adding the $(n-2)$- and $(n-1)$-handle to the RHB, $M_1$, of $B_1$ to produce $M$. But now $M_1 \hookrightarrow B_0$ induces a $\pi_1$ isomorphism, as $n-2$ and $n-1 \ge 3$.

\end{Proof}

However, $(B_1, M, M_1)$ is not the cobordism we seek.

\indexspace

Since the 1-handle $\alpha_1^1$ was trivially attached, we may attach a canceling 2-handle $\beta_1^2$. But now the original 2-handle $\alpha_1^2$ is trivially attached (by observation of its attaching loops with $t_1$ killed), so we may attach a canceling 3-handle $\beta_1^3$. Set $W_1^{n+1} = (M_1 \times \mathbb{I}) \cup \beta_1^2 \cup \beta_1^3$

\indexspace

But now, $B_0 \cup \beta_1^2 \cup \beta_1^3 = B_0\cup_{M_1} W_0 \approx M \times \mathbb{I}$. Invert $(W_0, M_1, M)$ so that it becomes $(W_0, M, M_1)$ with an $(n-3)$-handle $\gamma_1^{n-3}$ and an $(n-2)$-handle $\gamma_1^{n-2}$ attached to the RHB of $M \times \mathbb{I}$. Then $\iota_\#: \pi_1(M) \rightarrow \pi_1(W_0)$ is an isomorphism, as $n-3$ and $n-2 \ge 3$.

\indexspace

$(W_0, M, M_1)$ is the cobordism we seek.

\indexspace

(inductive step) \textit{Mutatis mutandis} the basis step. Start with $M_i^n$ and cross it with $\mathbb{I}$. Attach a trivial 1-handle $\alpha_i^1$ corresponding to $t_{i+1}$ to the right-hand boundary (RHB) of $M_i \times \mathbb{I}$. Attach a 2-handle $\alpha_i^2$ for the relator $t_{i+1} = [t_{i+1}, t_i]$ to the RHB of $M_i \times \mathbb{I}$. Set $B_i^{n+1}$ to be $(M_i \times \mathbb{I}) \cup \alpha_i^1 \cup \alpha_i^2$. 

\begin{Claim}
Then $M_{i+1} \hookrightarrow B_i$ induces a $\pi_1$ isomorphism.
\end{Claim}

\begin{Proof}

By inverting the handlebody decomposition, we may view $B_i$ as the result of adding the $(n-2)$- and $(n-1)$-handle to the RHB, $M_{i+1}$, of $B_i$ to produce $M_i$. But now $M_{i+1} \hookrightarrow B_i$ induces a $\pi_1$ isomorphism, as $n-2$ and $n-1 \ge 3$.

\end{Proof}

However, $(B_i, M_i, M_{i+1})$ is not the cobordism we seek.

\indexspace

Since the 1-handle $\alpha_i^1$ was trivially attached, we may attach a canceling 2-handle $\beta_i^2$. But now the original 2-handle $\alpha_i^2$ is trivially attached (by observation of its attaching loops with $t_{i+1}$ killed), so we may attach a canceling 3-handle $\beta_i^3$. Set $W_{i}^{n+1} = (M_{i+1} \times \mathbb{I}) \cup \beta_i^2 \cup \beta_i^3$

\indexspace

But now, $B_i \cup \beta_i^2 \cup \beta_i^3 = B_i \cup_{M_{i+1}} W_{i} \approx M_i \times \mathbb{I}$. Invert $(W_{i}, M_{i+1}, M_i)$ so that it becomes $(W_{i}, M_i, M_{i+1})$ with an $(n-3)$-handle $\gamma_i^{n-3}$ and an $(n-2)$-handle $\gamma_i^{n-2}$ attached to the RHB of $M_i \times \mathbb{I}$. Then $\pi_1(W_{i}) \cong G_i$, as as $n-3$ and $n-2 \ge 3$.

\indexspace

$(W_{i}, M_i, M_{i+1})$ is the cobordism we seek.

\indexspace

Set $W^{n+1} = W_0 \cup_{M_1} W_1 \cup_{M_2} W_2 \cup_{M_3} W_3 \cup_{M_5} \ldots$

\indexspace

Define neighborhoods of infinity $N_i = W_i \cup_{M_{i}} W_{i+1} \cup_{M_{i+1}} W_{i+2} \cup_{M_{i+2}} \ldots$ for $i \ge 1$. Then each $N_i$ is a connected neighborhood of infinity (clean except at $i = 1$).

(Pro-Fundamental Group System) Note that by properties b) and c) and a repeated application of Siefert-Van Kampen, the inverse sequence $\pi_1(N_0, p_0) \leftarrow \pi_1(N_1, p_1) \leftarrow \pi_1(n_2, p_2) \leftarrow \ldots$ is pro-isomorphic to $G_0 \leftarrow G_1 \leftarrow G_2 \leftarrow G_3 \leftarrow \ldots$. 

\indexspace

(Boundary) Also, clearly, $\partial W = M$.

\indexspace

(Each Clean Neighborhood has Finite Homotopy Type) 

It suffices to identify a cofinal sequence of clean neighborhoods of infinity having finite homotopy type. Toward that end, for each $i \ge 2$, let $N_i' = \beta_{i-1}^2 \cup N_i$, where $N_i = W_i \cup_{M_i} W_{i+1} \cup_{M_{i+1}} W_{i+2} \cup_{M_{i+2}} \ldots$.

\indexspace

The argument is complete when one shows that $\beta_{i-1}^2 \cup M_{i-1} \hookrightarrow N_i'$ is a homotopy equivalence.

\indexspace

This follows easily from the fact that for each $i \ge 1$, $W_i'$ strong deformation retracts onto $\beta_{i-1}^2 \cup M_i$, where $W'_i = W_i \cup \beta_{i-1}$. The proof of this fact is explained in \cite{G-T2}, Proposition 4.4.

\begin{Proof}

It suffices to show that $\beta_{i-1}^2 \cup M_i \hookrightarrow W_i'$ is a homotopy equivalence. Let $b_{i-1}^{n-2}$ be a belt disk for $\beta_{i-1}^2$ that intersects $\widetilde{\beta_{i-1}^2}$ in a belt disk $\widetilde{b_{i-1}^{n-2}}$, where $\widetilde{\beta_{i-1}^2}$ is a tubular neighborhood of the core disk of $\beta_{i-1}^2$ that sits inside $\hat{\beta_{i-1}^2}$.

\indexspace

By the Paint Can Lemma (which states that a cube less one face strong deformation retracts onto the remaining faces), $\beta_{i-1}^2 \cup M_i$ strong deformation retracts onto $b_{i-1}^{n-2} \cup M_i'$ ($\clubsuit$).

\indexspace

By a similar move, $W_i'$ strong denormalformation retracts onto $\widetilde{b_{i-1}^{n-2}} \cup W_i''$.

\indexspace

But $W_i''$ is a product, ($\ddag$), so we may collapse $\widetilde{b_{i-1}^{n-2}} \cup W_i''$ onto $b_{i-1}^{n-2} \cup M_i'$, which is a strong deformation retract of $\beta_{i-1}^2 \cup M_i$, by ($\clubsuit$).

\end{Proof}

To prove $N_i'$ strong deformation retracts onto $\beta_{i-1}^2 \cup M_i$ (that is, that $N_i'$ strong deformation retracts onto $\beta_{i-1}^2 \cup M_i$), collapse $W_i'$ onto $\beta_{i-1}^2 \cup M_i$, by the claim. But now, for $j > i$, $W_j'$ strong deformation retracts onto $(\beta_{j-1}^2 \cup M_j)$ extends (via the identity) to $(W_{j-1} \cup_{M_j} W_j)$ strong deformation retracts onto $W_{j-1}$, as $\beta_{j-1}^2 \subseteq W_j$. We may assemble all these strong deformation retractions to get a strong deformation retraction of $N_i'$ onto $\beta_{i-1}^2 \cup M_i$.

\end{Sketch}

%
%
%
%
%
%

  \singlespacing

\bibliography{dissertation} 
\bibliographystyle{plain}
\thispagestyle{references}

\end{document}